\documentclass[opre,nonblindrev]{informs3_arxiv}

\usepackage[backref]{enotez}
\let\footnote=\endnote
\let\enotesize=\normalsize
\def\notesname{Endnotes}%
\def\makeenmark{$^{\theenmark}$}
\def\enoteformat{\rightskip0pt\leftskip0pt\parindent=1.275em
  \leavevmode\llap{\makeenmark}}

\usepackage{natbib}
 \bibpunct[, ]{(}{)}{,}{a}{}{,}%
 \def\bibfont{\small}%
\TheoremsNumberedThrough     
\ECRepeatTheorems

\EquationsNumberedThrough    

\usepackage[utf8]{inputenc}
\usepackage[T1]{fontenc}
\usepackage{enumerate}
\usepackage{mathrsfs}

\usepackage{amsmath,natbib,amssymb}

\usepackage{booktabs} 

\usepackage{mathrsfs}
\usepackage{graphicx} 
\usepackage{bm}
\usepackage{dsfont}
\usepackage{fnpct}
\usepackage{mathtools}
\usepackage{comment}
\usepackage{mathrsfs}  
\mathtoolsset{showonlyrefs}
\usepackage{enumitem}
\usepackage{bbm}
\usepackage{makecell}

\usepackage{caption}
\usepackage{subcaption}

\usepackage{tocloft} 

\RequirePackage{hypernat}

\RequirePackage[colorlinks,citecolor=blue,urlcolor=blue]{hyperref}

\hypersetup{linktocpage}

\usepackage{bookmark}

\newtheorem{fact}{Fact}

\newcommand*{\QEDB}{\hfill\ensuremath{\square}}

\definecolor{DSgray}{cmyk}{0,0,0,0.7}
\definecolor{DSred}{cmyk}{0,0.7,0,0.7}
\definecolor{DSblue}{cmyk}{0.7, 0.7, 0, 0}

\begin{document}

\RUNAUTHOR{Hązła et al.}

\RUNTITLE{Bayesian Decision Making in Groups is Hard}

\TITLE{Bayesian Decision Making in Groups is Hard}



\ARTICLEAUTHORS{%
\AUTHOR{Jan Hązła$^{1,2}$, Ali Jadbabaie$^{1,3}$, Elchanan Mossel$^{1,2}$, M. Amin Rahimian$^{1}$}
\AFF{$^{1}$Institute for Data, Systems and Society \quad $^{2}$Department of Mathematics \quad $^{3}$Laboratory for Information and Decision Systems\\  Massachusetts Institute of Technology\\ \EMAIL{\{jhazla,jadbabai,elmos,rahimian\}@mit.edu}}
} 

\ABSTRACT{
  We study the computations that Bayesian agents undertake when exchanging opinions over a network. The agents act repeatedly on their private information and take myopic actions that maximize their expected utility according to a fully rational posterior belief. We show that such computations are NP-hard for two natural utility functions: one with binary actions, and another where agents reveal their posterior beliefs. In fact, we show that distinguishing between posteriors that are concentrated on different states of the world is NP-hard. Therefore, even approximating the Bayesian posterior beliefs is hard. We also describe a natural search algorithm to compute agents' actions, which we call \emph{elimination of impossible signals}, and show that if the network is transitive, the algorithm can be modified to run in polynomial time.}

\MSCCLASS{Primary: 91B06; secondary: 68Q25, 91A35, 62C10}

\ORMSCLASS{Primary: Games/group decisions: Voting/committees; secondary: Organizational studies: Decision making  Effectiveness/performance Information; Networks/graphs \\ {\hspace*{-\parindent}\emph{JEL}: D83, D85}. \\ {\hspace*{-\parindent}\emph{Working paper. \today. Authors are listed in alphabetic order.}}}

\KEYWORDS{Observational Learning, Bayesian Decision Theory, Computational Complexity, Group Decision-Making, Computational Social Choice, Inference over Graphs}

\maketitle

\section{Introduction}\label{sec:intro}

Many decision-making problems involve interactions among individuals
(agents) exchanging information with each other and striving to form rational opinions. Such situations arise in jury deliberations, expert committees, medical diagnoses, etc. Given their many important applications, relevant models of decision making in groups have been extensively considered over the years.

The first interesting case concerns only two agents. A fundamental insight offered by \cite{aumann1976agreeing} indicates that having common priors and a common knowledge of posterior beliefs imply agreement: Rational agents
cannot agree to disagree. Later work by \cite{geanakoplos1982we} demonstrates that such an agreement can be achieved in finite time, by broadcasting  posterior beliefs back and forth. \cite{Banerjee1992} and \cite{Bikhchandani1998} study the sequential interaction where each agent observes the decisions of everyone before her. \cite{acemoglu2011bayesian} extend the sequential learning model to a network environment where agents only observe actions of their neighbors (rather than all preceding actions). \cite{GaleANDkariv} consider repeated (rather than sequential) interactions over  social networks where agents update their beliefs after observing actions of each other. Following on, a large body of literature studies different aspects of rational opinion exchange, in particular the quality of information aggregation and learning in the limit
(cf. \cite{MT17,AO11} for two surveys of known results). 

Another prominent approach to the study of social learning  is to model non-Bayesian agents who use simpler, heuristic rules.
One reason for considering non-Bayesian heuristics (so called ``bounded rationality'') in place of fully rational updates is seeming intractability of Bayesian calculations: Information from different neighbors can exhibit complex correlations,
with no obvious way to account for, and remove them. For example, a Bayesian agent may have to account for the fact that her neighbors are influenced by the same source of information or even her own past actions \citep{eyster2014extensive,krishnamurthy2014online}.

Even though hardness of Bayesian computations in networked, learning models seems to be widely believed, we are not aware of any previous work making a rigorous argument for it. Our present work addresses this gap. We analyze algorithmic and complexity theoretic foundations of Bayesian social learning in two natural environments that are commonly studied in the literature. In one of them the actions broadcast by agents
are coarse, in the sense that they are single bits. In the other one, we assume that the actions are rich, consisting of agents'
full posterior beliefs. We show that the computations of the agents are intractable in both cases.

\subsection{Our contributions}

We analyze a fairly well-studied model of Bayesian social learning.
In this model there is a random variable $\theta$ which represents the unknown state of the world and determines payoffs from different actions. A network of agents receive private signals which are independent conditioned on the value of $\theta$. At every step $t = 0, 1, 2, \ldots$, each agent outputs an action $\mathbf{a}_{i,t}$ that maximizes her utility according to her current posterior distribution of $\theta$. The action is chosen myopically, i.e., only utility at the current time is considered and the posterior $\boldsymbol{\mu}_{i,t}$ is computed using Bayes rule. Agents learn actions of their neighbors on the network and proceed to the next step with updated posteriors.

For our hardness results, we study two natural variants of this model.
First, we consider the case of \emph{binary actions}, where the state, signals and actions are all binary, and each agent outputs the guess for the state $\theta\in\{0,1\}$ that is most likely according to her current belief.
This model can be thought of as repeated voting (e.g., during jury deliberations or the papal conclave in the Catholic Church). We are interested in the complexity of computations for producing Bayesian posterior beliefs  $\boldsymbol{\mu}_{i,t}$ or action $\mathbf{a}_{i,t}$. We also study the \emph{revealed belief model} where the utilities induce agents to reveal their current posteriors, or beliefs.

Following the detailed model description in Section~\ref{sec:BayesianModel}, we present our complexity results in Section~\ref{sec:NPhard}. We show that it is \emph{NP-hard} for the agents to compute their actions, both in the binary action and the revealed belief model. As a common tool in computational complexity theory, NP-hardness  provides rigorous evidence of \emph{worst-case} intractability.
Note that we only prove existence of  intractable network structures and private signals,
not that they are ``common'' or ``likely to arise''.
Also, our reductions critically rely on the network structure: They do not apply
to sequential models like the one in~\cite{Banerjee1992}.
 One might suspect that the beliefs can be efficiently approximated, even if they are difficult to compute exactly. This is unfortunately not the case, and we further prove a \emph{hardness-of-approximation} result: It is difficult even to distinguish between posterior beliefs that concentrate almost all of probability on one state and those that are concentrated on another state. In Section~\ref{sec:NPhard} we discuss in more detail what substantive economic assumptions are important in deriving our complexity results and some ways in which those results can be extended.

In Section \ref{sec:algorithms}, we study algorithms for Bayesian decision making in groups and describe a natural search algorithm to compute agents' actions. The Bayesian calculations are formalized as an algorithm for elimination of impossible signals (EIS), whereby the agent refines her knowledge by eliminating all profiles of private signals that are inconsistent with her observations. In Subsection \ref{sec:generalStrcutures}, we present recursive and iterative implementations of this algorithm. While the search over the possible signal profiles using this algorithm runs in exponential time, these calculations simplify in certain network structures. In Subsections \ref{sec:POSETs} and \ref{sec:algobelief}, we give examples of efficient algorithms for such cases. As a side result, we provide a partial answer to one of the questions raised by \cite{mossel2013making}, who provide an efficient algorithm for computing the Bayesian binary actions in a complete graph: We show that efficient computation is possible for other graphs that have a transitive structure when the action space is finite. In such \emph{transitive networks}, every neighbor of a neighbor of an agent is also her neighbor and therefore there are no indirect interactions to complicate the Bayesian inference.

\subsection{Related work}
\label{sec:literature}

Our results are related to the line of work that studies conditions for consensus and learning among rational agents \citep{smith2000pathological,mueller2013general, mossel2018social}. Consensus refers to all agents converging in their actions or belief (cf. \cite{GaleANDkariv,rosenberg2009informational} for consensus conditions in the network model that we study). Learning means that the consensus action is efficient, i.e., it represents the state of the world with high probability.
For example, \cite{MosselSlyTamuz14,mossel2015strategic} consider the binary action model (for myopic and forward-looking agents, respectively) and provide sufficient conditions for learning. These conditions are imposed on the network structure and consist of bounded out-degree and an ``egalitarian'' connectivity, whereby if an agent $i$ observes agent $j$, there is a reverse path from $j$ to $i$ of bounded length (this condition is trivially
satisfied for undirected networks).  

On the other hand,  positive computational results for Bayesian opinion exchange (including the analysis of short-run dynamics) are restricted to small networks (e.g., with three agents \cite[Section 5]{GaleANDkariv}, see also examples
in~\cite{rosenberg2009informational}) or special cases. The case of jointly Gaussian signals and beliefs exhibits a linear-algebraic structure that allows for tractable computations
(\cite{mossel2016efficient}, see also \cite{PersuasionBias}). \cite{dasaratha2018social} extend this setup to dynamic state spaces and
private signals. There are also efficient algorithms for special network structures, e.g., complete graphs and trees \citep{mossel2013making, TamuzTrees}.  Moreover, recursive techniques have been applied to analyze Bayesian decision problems with partial success, \cite{MosselSlyTamuz14,mossel2016efficient,harel2014more,TamuzTrees} and we also contribute to this literature by offering new cases where Bayesian decision making is tractable (cf. Subsections \ref{sec:POSETs} and \ref{sec:algobelief}).
This state of affairs might have to do with our computational hardness results. 

Other ways to achieve positive computational results are through alternative communication strategies or using non-Bayesian information exchange protocols. For example, \cite{acemoglu2014dynamics} analyze social learning among agents who directly communicate their entire information (represented as pairs of private signals and their sources). Since each piece of information is tagged, there is no confounding, and Bayesian updating is simple. On the other hand, the exchanged information has a significantly more complex form. In contrast, we think of our model as relevant to situations where, as is often the case, it is not practical to exhaustively list all of one's evidence and reasoning instead of stating or summarizing one's opinion. A popular approach to study bounded rationality is by replacing Bayesian actions with heuristic (non-Bayesian) rules \citep{DeGrootModel,BalaGoyal,GolubWisdomCrowd,li2017learning,molavi2018theory,Jadbabaie2012210,mueller2017general, arieli2019naive, arieli2019robust}. These rules are often rooted in empirically observed behavioral and cognitive biases. For example, \cite{li2017learning} consider a class of naive agents who take Bayesian actions but as if their local neighborhood is the entire network. This assumption removes the possibility of indirect interactions and, similar to the transitive structures (Subsection \ref{sec:POSETs}), simplifies Bayesian computations.  Our work is orthogonal and complementary to these studies. We prove that Bayesian reasoning is otherwise, in general, computationally intractable (because of the difficulty of delineating confounded sources of information). 

There are also works that focus on two agents estimating an arbitrary random variable \citep{aaronson2005complexity} --- this is in contrast to our model where the state of the world is correlated with the private signals in a simple way. The  computational result of \cite{aaronson2005complexity} concerns a protocol where the two agents keep exchanging their Bayesian posteriors with a deliberately added noise term. One might question how ``Bayesian'' such a protocol is, since the agents are not maximizing a utility function. On the other hand, the error terms can be reinterpreted as transmission noise or computation errors of rational agents (where the agents have common knowledge of the noise distribution). \cite{aaronson2005complexity} shows that this protocol can be efficiently implemented (approximately and on average with respect to private signals) for any constant number of rounds. As far as we can see, the proof of \cite{aaronson2005complexity} does not extend to many agents and networks. In Subsection~\ref{sec:noisy}, we show how to adapt our hardness reduction to this noisy action setting. Notwithstanding, we cannot logically exclude the possibility of a result like \cite{aaronson2005complexity}, since we show only worst-case hardness and the algorithm in \cite{aaronson2005complexity} works on average. Therefore, we leave it as an interesting open problem: In the network model with noise, does there exist an average-case efficient algorithm, or are computations hard on average (at least with respect to private signal profiles)?

In fact, our results can be also interpreted in the context of other works pointing at computational reasons for why economic or sociological models fail to accurately reflect reality (cf., e.g., \cite{ABBG} on the computational complexity of financial derivatives and \cite{velupillai2000computable} on the computable foundations of economics). On the one hand, a model cannot be considered plausible if it requires the participants or agents to perform computations that need a prohibitively long time. On the other hand, the predictions of such a model can be rendered inaccessible by the computational barriers. The literature on computational hardness of Bayesian reasoning in social networks is nascent. There are some hardness results in the literature on Bayesian inference in graphical models (see \cite{Kwisthout11} and references therein), but these are quite different from models considered in this work. \cite{papadimitriou1987complexity} consider partially observed Markov decision processes (POMDP). These Markovian processes are not directly comparable to our model, but they exhibit similar flavor in so far as repeated interactions are concerned. \cite{papadimitriou1987complexity} prove that computing optimal expected utility in a POMDP is PSPACE-hard, achieving a stronger notion of hardness than NP-hardness. However, their result does not extend to hardness of approximation, i.e., they only show that it is hard to decide if the optimal agent's strategy achieves positive (but possibly very small) expected utility. Moreover, the setup for Bayesian decision making in groups is different
(arguably less general, i.e., more challenging for a hardness proof) than a POMDP. Subsequently, we need different techniques for our purposes.

We also point out a follow-up work by the authors of this paper \citep{BayesPSPACE18}, where we use significantly more technical arguments to show that the computations in the binary action model are also (worst case) PSPACE-hard to approximate. We believe the details of the latter work might be of interest to complexity theorists. Here, we focus on developing more general arguments to inform operations research and social learning applications.

\section{The Bayesian Group Decision Model}\label{sec:BayesianModel}

We consider a finite group of agents, whose interactions are represented by a fixed directed graph $\mathcal{G}$.
For each agent $i$ in $\mathcal{G}$, $\mathcal{N}_i$ denotes her neighborhood:
The subset of agents whose actions are observed by agent $i$.
Without loss of generality, we will assume that $i\in\mathcal{N}_i$, i.e.,
an agent always observes herself.

We model the topic of the discussion/group decision process by a state $\theta$ belonging to a finite set $\Theta$. For example, in the course of a political debate, $\Theta$ can be the set of all political parties with $\theta$ representing
the party that is most likely to increase
society's welfare.
The value of $\theta$ is not known to the agents, but
they all start with a common prior belief about it,
which is a distribution with probability mass function
$\nu(\mathord{\cdot}):\Theta \to [0,1]$.

Initially, each agent $i$ receives a private signal
$\mathbf{s}_i$, correlated with the state $\theta$.
The private signal  $\mathbf{s}_i $ belongs to a finite set $\mathcal{S}_i$
and its distribution conditioned on $\theta$ is denoted by ${\mathbb{P}}_{i,\theta}(\mathord{\cdot})$,
which is referred to as the \emph{signal structure} of agent $i$.
Conditioned on the state $\theta$, the signals $\mathbf{s}_i$ are independent
across agents, and
we use ${\mathbb{P}}_{\theta}(\mathord{\cdot})=\prod_{i}\mathbb{P}_{i,\theta}(\cdot)$ to denote their joint product
distribution.

After receiving the signals, the agents
interact repeatedly,  in discrete times $t=0,1,2,\ldots$. Associated with every agent $i$ is an action space $\mathcal{A}_i$ that represents the choices available to her at any time $t\in\mathbb{N}_0$, and a utility $u_i(\mathord{\cdot},\mathord{\cdot}):\mathcal{A}_i\times\Theta \to \mathbb{R}$ which represents her preferences with respect to combinations of actions and
states. At every time $t\in\mathbb{N}$, agent $i$ takes action $\mathbf{a}_{i,t}$
that maximizes  her expected utility based on her \emph{observation history}
$\mathbf{h}_{i,t}$:
\begin{align}
  \mathbf{a}_{i,t} = \argmax_{a_i\in\mathcal{A}_i}
\mathbb{E}[u_{i}(a_i,\theta)\mid\mathbf{h_{i,t}}],
  \label{eq:01}
\end{align}
where the history $\mathbf{h}_{i,t}$ is defined as
$\{\mathbf{s}_i\} \cup \{\mathbf{a}_{j,\tau}$
 for all $j\in\mathcal{N}_i$, and $\tau < t\}$,
i.e., agent $i$ observes her
private signal, as well as actions of all her neighbors at times strictly less
than $t$.

The network,
signal structures, action spaces and utilities, as well as the prior,
are all common knowledge among the agents.
We use the notation $\argmax_{a \in \mathcal{A}}$ to include the following, common knowledge, rule when the maximizer is not unique: We assume that the action spaces are (arbitrarily) ordered and an agent will break ties by choosing the lowest-ranked action in her ordering.
The specific tie-breaking rule is not important for our results.
The agents' behavior is myopic 
in that it does not take into account
strategic considerations about future rounds; cf. Subsection \ref{sec:forward-looking}.

We denote the Bayesian posterior belief of agent $i$ given her history of observations by its probability mass function $\boldsymbol\mu_{i,t}(\mathord{\cdot}):\Theta \to [0,1]$.
In this notation, the expectation in \eqref{eq:01} is taken with respect to the Bayesian posterior belief $\boldsymbol{\mu}_{i,t}$.

\begin{figure}[t]
\centering
\includegraphics[scale=0.5]{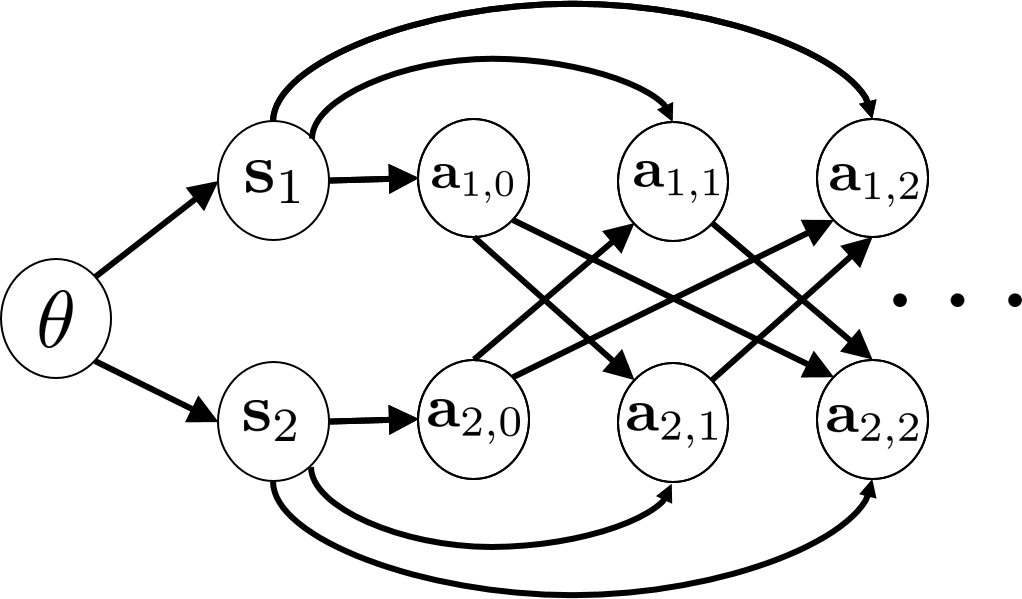}
\caption{The Decision Flow Diagram for Two Bayesian Agents}
\label{fig:DFD}
\end{figure}

To sum up, agent $i$ at time $t$ chooses an action $\mathbf{a}_{i,t} \in \mathcal{A}_i$,
maximizing her expected utility conditioned on the observation history $\mathbf{h}_{i,t}$.
Then, she observes the most recent actions of her neighbors
$\{\mathbf{a}_{j,t}$ for all $j\in\mathcal{N}_i\}$, updates her action to $\mathbf{a}_{i,t+1} \in \mathcal{A}_i$, and so on. A decision flow diagram for an example of two interacting agents is
provided in Figure~\ref{fig:DFD}.

Our main focus in this paper is  on the computational and algorithmic aspects
of the group decision process. Specifically, we will be concerned with the
following computational problem:
\begin{problem}[GROUP-DECISION] At a time $t$, given the graph structure $\mathcal{G}$, agent $i$ and the observation history $\mathbf{h}_{i,t}$, determine the Bayesian action ${\mathbf{a}}_{i,t}$.
\end{problem}

\subsection{Natural Utility Functions: Binary Actions and Revealed Beliefs}\label{sec:revealedBeliefs}

A natural example of a utility function is based on the idea of repeated
voting, for example, as an idealized model of jury deliberations or the papal conclave in the Catholic Church. In this model, the possible actions correspond to the states of the world, i.e.,
$\mathcal{A}_i = \Theta$ and
the utilities are given by
$u_i(a, \theta) = \mathbbm{1}(a = \theta)$.
In other words, the agents receive a unit reward for guessing the state correctly and zero otherwise.
The expected reward of agent $i$ at time $t$ is
maximized by choosing the action that corresponds to the maximum probability
in ${\boldsymbol\mu}_{i,t}$, i.e. the maximum a posteriori probability (MAP) estimate. In case of binary world $\Theta = \{0,1\}$ with uniform prior and binary private signals
$\mathcal{S}_i=\{0,1\}$ we call this example the \emph{binary action} model. 

In another important example, which we call the \emph{revealed belief} model, the agents reveal their complete posteriors,
i.e., $\boldsymbol\mu_{i,t}$.
Formally, let
$\Theta:=\{\theta_1,\ldots,\theta_{m}\}$ and let
$\overline{e}_j \in \mathbb{R}^{m}$ be a column vector of all zeros except for
its $j$-th element which is equal to one. Furthermore, we relax the requirement that the action spaces $\mathcal{A}_i$ are finite sets; instead, for each agent $i\in [n]$ let $\mathcal{A}_i$ be the $m$-dimensional probability simplex: $\mathcal{A}_i = \{(x_1,\ldots,x_m)^T \in \mathbb{R}^m:\sum_{i=1}^{m}x_i = 1 \, \mbox{ and } \, x_i \geq 0, \forall i \}$. If the utility assigned to an action $\overline{a}:=(a_1,\ldots,a_m)^T\in\mathcal{A}_i$ and a state $\theta_j\in\Theta$ measures the squared Euclidean distance between $\overline{a}$ and $\overline{e}_j$, then it is optimal for agent $i$ to reveal her belief  $\mathbf{a}_{i,t} = ({\boldsymbol\mu}_{i,t}({\theta}_1),\ldots,{\boldsymbol\mu}_{i,t}({\theta}_m))^T$. We can state a special case of the GROUP-DECISION model in the revealed belief
setting:

\begin{problem}[GROUP-DECISION with revealed beliefs] At any time $t$, given the graph structure $\mathcal{G}$, agent $i$ and the observation history
  $\mathbf{h}_{i,t}$,
determine the Bayesian posterior belief ${{\boldsymbol\mu}}_{i,t}$.
\end{problem}

\subsection{Log-Likelihood Ratio and Log-Belief Ratio Notations} \label{sec:loglike}

Consider a finite state space $\Theta = \{\theta_1,\ldots,\theta_m\}$ and for all $2\leq k \leq m$ and $s\in\mathcal{S}_i$, let:

\begin{align}\label{eq:11}
 {\lambda}_i(s, \theta_k) := \log\left(\frac{{\mathbb{P}}_{i,\theta_k}(s)}{{\mathbb{P}}_{i,\theta_1}(s)}\right) \, , \,  \boldsymbol{\phi}_{i,t}(\theta_k) := \log\left(\frac{\boldsymbol{\mu}_{i,t}(\theta_k)}{\boldsymbol{\mu}_{i,t}(\theta_1)}\right) \, , \,
 {\gamma}(\theta_k) := \log\left(\frac{{\nu}(\theta_k)}{{\nu}(\theta_1)}\right).
\end{align}
We will also write
$\boldsymbol{\lambda}_i(\theta_k):={\lambda}_i(\mathbf{s}_i,\theta_k)$.
We will call $\boldsymbol\lambda_i$ the \emph{(signal) $\log$-likelihood ratio}
and $\boldsymbol\phi_{i,t}$ the \emph{$\log$-belief ratio}.
If we assume that the agents start from uniform prior beliefs and the size of
the state space is $m=2$ (as will be the case for the hardness results
in Section~\ref{sec:NPhard}), we can employ a simpler notation. First, with uniform priors, we have $ {\gamma}(\theta_k) = \log\left({{\nu}(\theta_k)}/{{\nu}(\theta_1)}\right) = 0$ for all $k$. Moreover, with binary state space $\Theta = \{0, 1\}$ we only need to keep track of one set of $\log$-belief and $\log$-likelihood ratios $\boldsymbol{\lambda}_i := \boldsymbol{\lambda}_i(1) = \log\left({{\mathbb{P}}_{i,1}(\mathbf{s}_i)}/{{\mathbb{P}}_{i,0}(\mathbf{s}_i)}\right)$, and $\boldsymbol{\phi}_{i,t} = \boldsymbol{\phi}_{i,t}(1) =   \log\left({\boldsymbol{\mu}_{i,t}(1)}/{\boldsymbol{\mu}_{i,t}(0)}\right)$. Henceforth, we use $\boldsymbol{\lambda}_i$ and $\boldsymbol{\phi}_{i,t}$ as there is no risk of confusion in dropping their arguments. 

Note that in the setting with binary state and signals
($\mathcal{S}_i=\{0,1\}$), there is a
one-to-one correspondence between informative signal structures satisfying  $\mathbb{P}_{i,0}(1)\ne\mathbb{P}_{i,1}(1)$, and $\log$-likelihood ratios satisfying $\boldsymbol{\lambda}_{i}(0)\cdot\boldsymbol{\lambda}_{i}(1)<0$.
Accordingly, we sometimes use $\log$-likelihood ratios to specify signal structures.

\begin{example}[Belief Exchange in the First Two Rounds]\label{example:initial_rounds}
 To give some intuition about our model and illustrate the usefulness of the $\log$-likelihood ratio and $\log$-belief ratio 
notations, we explain how the agents in the binary action model can compute their actions at $t=0$ and $t=1$.
We consider informative binary private signals $\mathbf{s}_i \in \{0, 1\}$ with $\mathbb{P}_{i,1}(1) > \mathbb{P}_{i, 0}(1)$. We focus on computing the $\log$-likelihood ratio ($\boldsymbol\phi_{i,t}$), since $\mathbf{a}_{i,t} = 1$ if, and only if,
$\boldsymbol\phi_{i,t} > 0$.

At time zero, the posterior and $\log$-belief ratio of agent $i$ are determined by her private signal, as follows:
\begin{align}
  \boldsymbol\mu_{i,0}(1) =
  \frac{\mathbb{P}_{i,1}(\mathbf{s}_i)}{\mathbb{P}_{i,0}(\mathbf{s}_i)+
  \mathbb{P}_{i,1}(\mathbf{s}_i)}\; ,\qquad\qquad
  \boldsymbol\phi_{i,0} = \log \left(
  \frac{\mathbb{P}_{i,1}(\mathbf{s}_i)}{\mathbb{P}_{i,0}(\mathbf{s}_i)}
  \right)\;.
\end{align}
Therefore, we get $\mathbf{a}_{i,0} = \mathbf{s}_i$ since $\mathbb{P}_{i,1}(1) > \mathbb{P}_{i, 0}(1)$. At time one, agent $i$ observes the actions, and therefore infers the private signals, of her neighbors. Since the private signals are conditionally
independent, the respective $\log$-likelihood ratios add up and we get the following
expression (recall that $i\in\mathcal{N}_i$):
\begin{align}
  \boldsymbol\phi_{i,1} =
  \sum_{j\in\mathcal{N}_i} \boldsymbol\phi_{j,0}=
  \sum_{j \in \mathcal{N}_i} \log \left(
  \frac{\mathbb{P}_{j,1}(\mathbf{a}_{j,0})}{\mathbb{P}_{j,0}(\mathbf{a}_{j,0})}
  \right) =  \sum_{j \in \mathcal{N}_i} \boldsymbol\lambda_{j} \; .
  \label{eq:beliefexchange}
\end{align} However, the computation becomes significantly more involved at later times. This is because one needs to account for dependencies and redundancies in agents' information and the resulting actions.
\end{example}

\section{Hardness of Bayesian Decisions}\label{sec:NPhard}

Our hardness results use a standard approach from complexity theory; cf., e.g.,~\cite{AB09}.
We establish \emph{NP-hardness} of computations in both binary action and
revealed belief models. We do so by exhibiting reductions from problems
that are known to be NP-hard. As shown below, two covering problems:
vertex cover and set cover, turn out to be convenient starting
points for our reductions. We now present our main hardness results.

{\begin{figure}[t]
\centering
\begin{subfigure}[b]{0.3\textwidth}
   \includegraphics[width=1\linewidth]{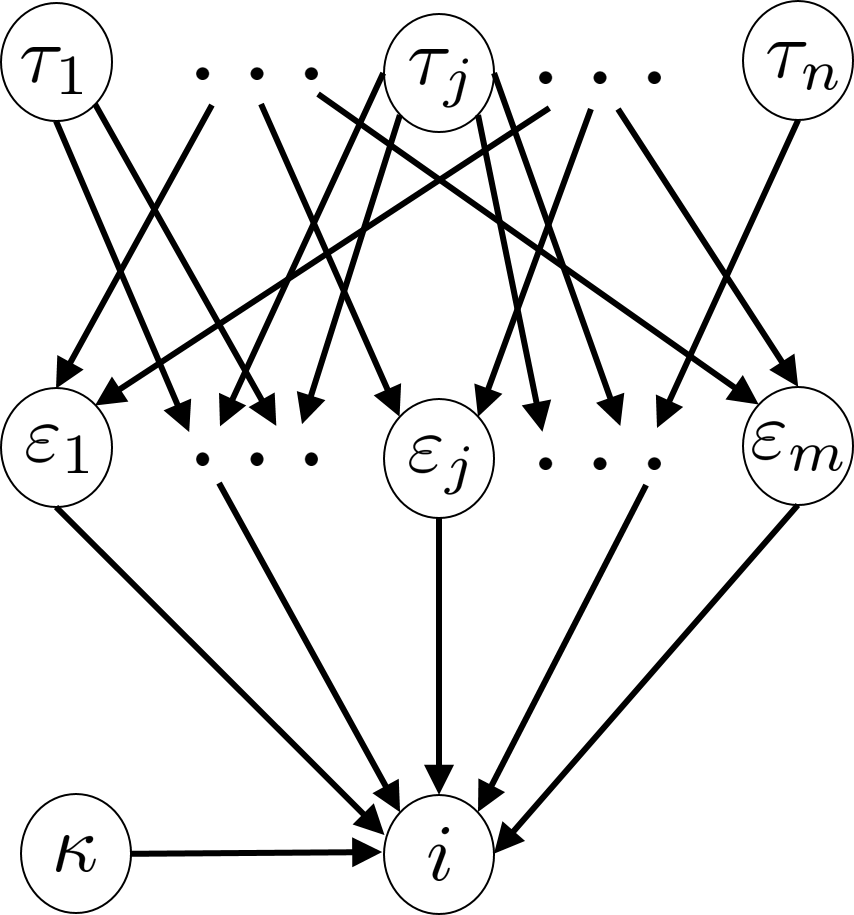}
   \caption{ }
   \label{fig:VERTEXCOVER}
\end{subfigure}~
\begin{subfigure}[b]{0.3\textwidth}
  \includegraphics[width=1\linewidth]{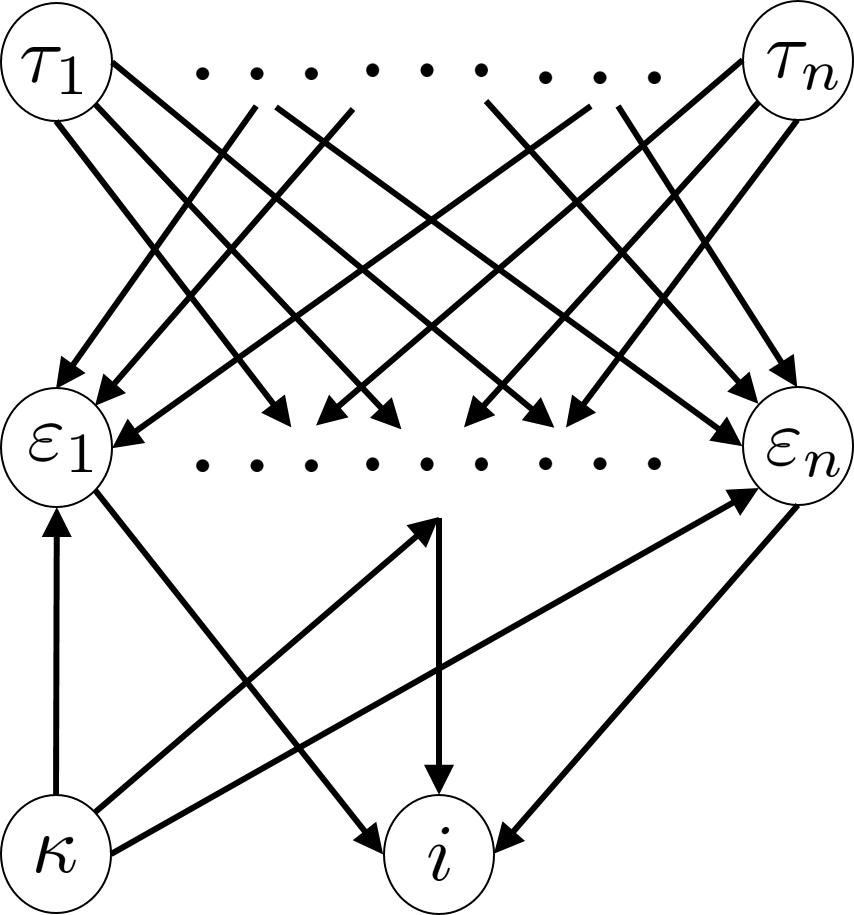}
\caption{ }
\label{fig:EXACTCOVER}
\end{subfigure}
\caption{(\subref{fig:VERTEXCOVER}) Illustration of the VERTEX-COVER reduction (Theorem \ref{THM:NPHARDNESSACTIONS}); every edge $\varepsilon_j$ is connected to its two vertices, and every vertex is connected to all its incident edges. (\subref{fig:EXACTCOVER}) Illustration of the EXACT-COVER reduction (Theorem \ref{THM:NPHARDNESSBELIEFS}); every element $\varepsilon_j$ belongs to exactly three sets and every set $\tau_j$ contains exactly three elements.}
\end{figure}
} 

\begin{theorem}[Binary Action Model]\label{THM:NPHARDNESSACTIONS}
  The GROUP-DECISION problem in the binary action model is NP-hard
  at $t=2$. Furthermore, for a network of $n$ Bayesian agents in the binary action model, it is NP-hard to distinguish between
  posterior beliefs $\boldsymbol{\mu}_{i,2}(0)<\exp(-\Omega(n))$
  and $\boldsymbol{\mu}_{i,2}(1)<\exp(-\Omega(n))$.
\end{theorem}

\proof{Proof sketch of Theorem~\ref{THM:NPHARDNESSACTIONS}.} Appendix \ref{App:vertexConverProof} contains a detailed proof. Our reduction is from an NP-hard problem of \emph{approximating vertex cover} (VERTEX-COVER). A vertex cover on an undirected graph $\hat{\mathcal{G}}_{m,n}$ with $n$ vertices and $m$ edges is a subset of vertices (denoted by $\hat{\Sigma}$), such that each edge touches at least one vertex in $\hat{\Sigma}$. We consider the approximation version of VERTEX-COVER, where every input graph belongs to one of two cases: 
\begin{enumerate}[label=(\roman*)]
\item the YES case, where it has at least one small vertex cover (say, smaller than $0.85n$),
\item the NO case, where all its vertex covers are large (say, larger than $0.999n$).
\end{enumerate} It is NP-hard to distinguish between these two cases.

We show an efficient reduction that maps a graph $\hat{\mathcal{G}}_{m,n}$ to an instance of GROUP-DECISION in the binary action model.
We encode the structure of $\hat{\mathcal{G}}_{m,n}$ by a two-layer network, where the first layer is comprised of ``vertex agents'', which are connected to
``edge agents'' in the second layer based on the incidence relations in $\hat{\mathcal{G}}_{m,n}$ (see Figure~\ref{fig:VERTEXCOVER}).

We let the vertex agents $\tau_1,\ldots,\tau_n$ receive Bernoulli private signals
with signal structure given by $\overline{p}:=\mathbb{P}_{\tau_i,1}(1)=0.4$ and $\underline{p}:=\mathbb{P}_{\tau_i,0}(1)=0.3$. Each edge agent
($\varepsilon_j$) observes two vertex agents corresponding to its incident vertices in $\hat{\mathcal{G}}_{m,n}$. The private signals of edge agents are uninformative. We can verify that since $\overline{p}(1-\overline{p}) = 0.24 > 0.21 = \underline{p}(1-\underline{p})$,
an edge agent $\varepsilon_{j}$ takes action one at time one  ($\mathbf{a}_{\varepsilon_{j},1} = 1$) if, and only if, at least one of the two neighboring vertex agents ($\tau_i$) receives private signal $\mathbf{s}_{\tau_i} = 1$.

Agent $i$ (whose decision we show to be NP-hard) receives an uninformative private signal, and observes all edge agents as well (see Figure~\ref{fig:VERTEXCOVER}). To complete the reduction, we need
to specify the observation history of agent $i$, and we do so by saying
that all edge agents
announce action one at time one $\mathbf{a}_{\varepsilon_j,1} = 1$.
By our previous observation, this is equivalent to saying that
the private signals of vertex agents form a vertex cover of $\hat{\mathcal{G}}_{m,n}$.

The crux of the proof is in showing the following property:
\begin{itemize}
\item If every vertex cover of
  $\hat{\mathcal{G}}_{m,n}$ has size at least $0.999n$,
  then agent $i$ concludes that at least $0.999n$ of vertex agent private
  signals are ones.
\item On the other hand,
  if $\hat{\mathcal{G}}_{m,n}$ has a vertex cover of size at most $0.85n$,
  then agent $i$ concludes that, almost certainly, at most $0.998n$
  of private signals are ones.
\end{itemize}

The first statement is clear.
However, if there exists a vertex cover of size $0.85n$, the private
signals might come from this small vertex cover just as well as from any of
the larger covers. Since
$\underline{p}=0.3$ and $\overline{p}=0.4$, the size of any vertex cover
is much larger than expected number of ones among
the private signals,
regardless of the state $\theta$.
One could hope that the concentration of measure would imply that seeing a smaller vertex cover is relatively much more likely, even
if there is a significantly greater total number of large vertex covers.
In Appendix~\ref{App:vertexConverProof}, we use a Chernoff bound to conclude that this is indeed the case, and agent $i$ can infer
that, almost certainly, the private signals form a vertex cover of size at
most $0.998n$.

After establishing that it is NP-hard to distinguish between at least $0.999n$ ones and at most $0.998n$ ones among the private signals, our construction concludes with a simple trick. We will explain the idea assuming
a gap between $0.8n$ and $0.6n$ instead of inconveniently small
$0.999n$ and $0.998n$. The complete details are provided in Appendix~\ref{App:vertexConverProof}.

Assume that agent $i$ additionally observes another agent $\kappa$.
Agent $\kappa$ does not observe anyone and reveals to agent $i$
a very strong, independent private signal equivalent to $n$ signals
of vertex
agents, all of them with value zero.
If agent $i$ is in the case where at least $0.8n$ vertex signals are ones,
then her total observed signal strength is equal to at least
$0.8n$ ones out of $2n$ total, i.e., at least $40\%$ of all signals
are ones. Given that $\overline{p}=0.4 = 40\%$, agent $i$ concludes that
almost certainly $\theta=1$, i.e., $\boldsymbol{\mu}_{i,2}(0)\approx 0$.
On the other hand, in case where (almost certainly)
at most $0.6n$ vertex signals are ones, total signal strength
is at most $0.6n$ out of $2n$, i.e., $30\%$ of possible signals 
and, recalling $\underline{p}=0.3=30\%$, agent $i$ concludes that
$\boldsymbol{\mu}_{i,2}(1)\approx 0$.\endnote{Technically we showed coNP-hardness, i.e., our reduction mapped instances with small vertex cover onto GROUP-DECISION instances with $\theta = 0$ and instances with only large vertex covers onto GROUP-DECISION with $\theta=1$. However, due to the symmetric nature of GROUP-DECISION, NP-hardness is immediately obtained by inverting the meanings of $0$ and $1$ labels of states and private signals.  In particular, since GROUP-DECISION at $t=2$ is both NP-hard and coNP-hard, it  is likely to be strictly harder than NP-complete (see~\cite{AB09}).} \QEDB
\endproof

\begin{remark}
  A priori one might suspect that the difficulty of distinguishing between
  $\mathbf{a}_{i,2} = 0$ and $\mathbf{a}_{i,2} = 1$ arises only if the belief
  of agent $i$ is very close to the threshold $\boldsymbol\mu_{i,2} \approx 1/2$.
  However, in our reduction the opposite
  is true: For a computationally bounded agent, it is hopeless to distinguish
  between worlds where $\theta = 0$ with high probability (w.h.p.), and $\theta = 1$
  w.h.p. This can be thought of as a strong hardness of approximation result.
\end{remark}

We also have a matching result for the revealed belief model:

\begin{theorem}[Approximating Beliefs]\label{THM:NPHARDNESSBELIEFS} The GROUP-DECISION problem is NP-hard in the revealed belief model  with uniform priors,  binary states $\Theta = \{0,1\}$, and binary private signals $\mathbf{s}_i \in \{0, 1\}$. In particular, for a network of $n$ Bayesian agents at $t=2$, it is NP-hard to distinguish between beliefs $\boldsymbol\mu_{i,2}(0) \le \exp(-\Omega(n))$ and $\boldsymbol\mu_{i,2}(1) \le \exp(-\Omega(n))$.
\end{theorem}

\proof{Proof sketch for Theorem \ref{THM:NPHARDNESSBELIEFS}.} Appendix \ref{App:exactConverProof} contains a detailed proof. Our reduction is from a variant of an NP-complete problem EXACT-COVER. Let $n$ be a multiple of three
and consider a set of $n$ elements $\hat{\mathcal{E}}_n = \{\varepsilon_1,\ldots,\varepsilon_n\}$ and a family of $n$ subsets of $\hat{\mathcal{E}}_n$ denoted by $\hat{\mathcal{T}}_n = \{\tau_1,\ldots,\tau_n\}$, ${\tau}_j\subset\hat{\mathcal{E}}_n$ for all $j\in[n]$. EXACT-COVER is the problem of deciding if there exists a collection $\hat{\mathcal{T}} \subseteq \hat{\mathcal{T}}_n$ that exactly covers $\hat{\mathcal{E}}_n$, that is, each element $\varepsilon_i$ belongs to exactly one set in $\hat{\mathcal{T}}$. We use a restriction of EXACT-COVER where each set has size three and each element appears in exactly three sets; hence, if the exact cover exists,
then it consist of $n/3$ sets.
  
We use a two-layer network to encode the inclusion relations between the elements $\hat{\mathcal{E}}_n$ and subsets $\hat{\mathcal{T}}_n$. There are $n$ agents
$\tau_1,\ldots,\tau_n$ in the first layer to encode the subsets and $n$ agents
$\varepsilon_1,\ldots,\varepsilon_n$
in the second layer to encode the elements. Each ``element agent'' observes
three ``subset agents'' corresponding to subsets to which the element belongs (see Figure \ref{fig:EXACTCOVER}). Agent $i$ (whose decision we show to be NP-hard) observes the reported beliefs of all element agents. There is also one auxiliary agent $\kappa$ that is observed by all element agents.

The private signals of agent $i$ and the element agents are non-informative. The subset agents observe i.i.d.~binary signals and the auxiliary agent $\kappa$ observes another independent binary signal, but with a different distribution. We set up the signal structures
and the beliefs transmitted by the element agents to agent $i$ such that there are two possible outcomes: Either $\mathbf{s}_\kappa = 0$ and all subset agents received positive signals $\mathbf{s}_{\tau_i} = 1$; or, $\mathbf{s}_{\kappa} = 1$ and the private signals of subset agents form an exact cover of the elements. Of course, the second alternative is possible only if an exact cover exists.

The first alternative implies that all subset agents received ones as private signals,
and therefore $\theta = 1$ with high probability.
In case of the second alternative, we show that almost certainly
only one-third of subset agents received ones, and therefore
$\theta = 0$ with high probability. Therefore,
if there is no exact cover, agent $i$ should compute $\boldsymbol\mu_{i,2}(0) \approx 0$ and otherwise $\boldsymbol\mu_{i,2}(1) \approx 0$. \QEDB

\endproof

We conclude this section by discussing some aspects and limitations of our
proof. We also examine the economic assumptions behind our results
and discuss what happens when these assumptions are relaxed.

\subsection{Worst-Case and Average-Case Reductions}
\label{sec:average}

Our reductions are worst-case, both with respect to networks and signal
profiles. That is, we show hardness only for a specific class of networks,
and for signal profiles in those networks that arise with exponentially small
probability. We cannot exclude existence of an efficient algorithm that computes
Bayesian beliefs for all network structures, with high probability
over signal profiles. Notwithstanding, any such purported algorithm must have
a good reason to fail on our hard instances.

This reflects a general phenomenon in  computational complexity, where average-case hardness, even when suspected to hold,
seems to be significantly more difficult to rigorously demonstrate
(see \cite{bogdanov2006average} for one survey). We leave as a fascinating open
problem if our results can be improved, for example for worst-case
networks and average-case signal profiles. One thing to note in this regard
is that our reductions encode the witnesses to NP problems
(vertex and set covers) as signal profiles. That necessarily means that
for hard positive instances (e.g., graphs with a small vertex cover) relevant
signal profiles will arise only with tiny probability: Otherwise these
instances would be easy to solve by sampling a potential witness at random. Significant new ideas might be needed to overcome this problem.

On the positive side, the worst-case nature of our hard instances makes it
potentially easier to embed them in more general or modified settings. We discuss several concrete cases below.

\subsection{Forward-looking Agents}\label{sec:forward-looking}
Our results are restricted to myopic agents. In the general framework of forward-looking utility maximizers with discount factor $\delta$, myopic agents are obtained as
a special case by completely discounting the future pay-offs ($\delta \to 0$).

The computational difficulties for strategic agents seem to be at least as large
as for myopic agents, but we do not offer any formal results.
Due to the multiplicity of equilibria suggested by the folk theorem
\citep{fudenberg1986folk} --- see also examples in
\cite{rosenberg2009informational} and \cite{mossel2015strategic} --- it is unclear
to us how to make the computational problem well-posed. On the other hand,
since in the limit $t\to\infty$ the agents in any equilibrium act myopically
\citep{rosenberg2009informational}, it seems plausible to expect that
their computations will be similarly hard as in our analysis.


\subsection{Directed Links} In both our reductions, we use directed acyclic
graphs. This is arguably a simpler case from an inference viewpoint, since in
networks containing cycles (including those with bidirectional links) an agent needs to take into account her own, possibly indirect,
influence on her neighbors. Therefore, our hardness results hold true, in spite of the (simpler) acyclic structure of our hard examples.

In effect, our hardness results are applicable to undirected (bidirectional) networks
without loss of generality. The reason is that replacing directed links with bidirectional ones does not affect any relevant inferences in our reductions. In particular, our results apply to networks that exhibit agreement and learning,
cf.~\cite{MosselSlyTamuz14}. It is worth noting that since our results are achieved in a basic model with binary state and private signals, they can be easily embedded
in richer settings, e.g., with signal structures given by continuous
distributions.

{\begin{figure}[t]
\centering
\begin{subfigure}[b]{0.17\textwidth}
   \includegraphics[width=1\linewidth]{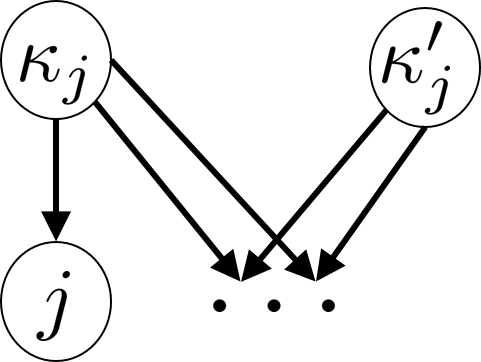}
   \caption{ }
   \label{fig:auxillary_non_uniform_priors}
\end{subfigure}~
\begin{subfigure}[b]{0.25\textwidth}
   \includegraphics[width=1\linewidth]{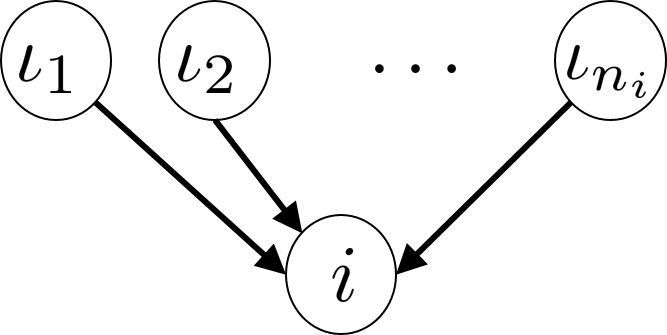}
   \caption{ }
   \label{fig:auxillary_n_i}
\end{subfigure}~
\begin{subfigure}[b]{0.29\textwidth}
  \includegraphics[width=1\linewidth]{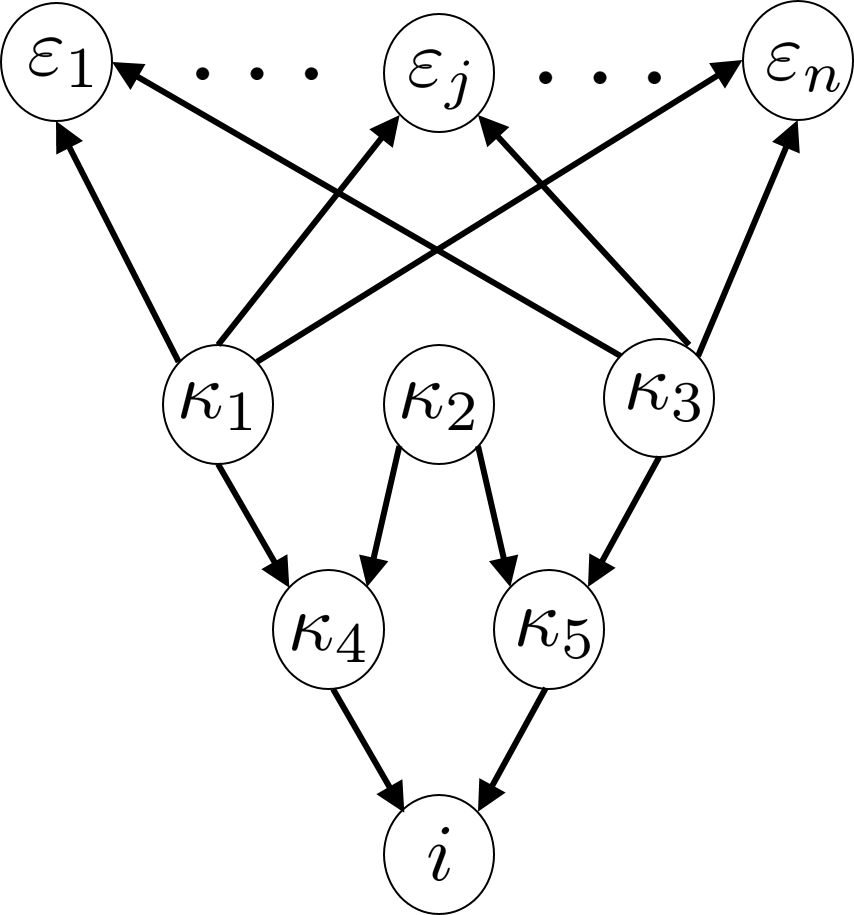}
\caption{ }
\label{fig:auxillary_iid_rx3}
\end{subfigure}
\begin{subfigure}[b]{0.18\textwidth}
  \includegraphics[width=1\linewidth]{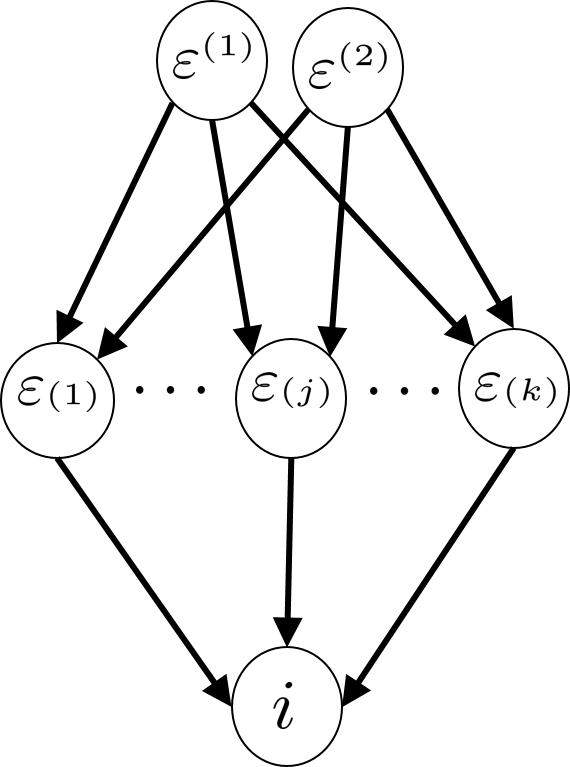}
\caption{ }
\label{fig:auxillary_noisy}
\end{subfigure}
\caption{(\subref{fig:auxillary_non_uniform_priors}) We can cancel out the effect of a distinct, non-uniform prior in an agent $j$ by adding two auxiliary agents ($\kappa_j$ and $\kappa'_j$), and have agent $j$ observe only one of them. Both added agents will be observed by every other agent. (\subref{fig:auxillary_n_i}) We can replace the auxiliary agent $\kappa$ in the VERTEX-COVER reduction by $n_i$ agents $\iota_1$, $\ldots$, $\iota_{n_i}$ with zero signals drawn from the same i.i.d distribution as the vertex agents. (\subref{fig:auxillary_iid_rx3}) We can replace the auxiliary agent $\kappa$ in the EXACT-COVER reduction by five agents $\kappa_1$, $\ldots$, $\kappa_5$ with i.i.d. signals and set up their received signals and the observation structure such that the signals of $\kappa_1$ and $\kappa_3$ necessarily agree. (\subref{fig:auxillary_noisy}) We can modify the VERTEX-COVER reduction to work with noisy binary actions. Here each pair of vertex agents ($\varepsilon^{(1)}$ and $\varepsilon^{(2)}$) are observed by a collection of edge agents $\varepsilon{(1)}$, $\dots$, $\varepsilon{(k)}$ who all report the same noisy actions $\mathbf{a}'_{\varepsilon{(1)},1} = \ldots = \mathbf{a}'_{\varepsilon{(k)},1} = 1$.}
\label{fig:auxillary_extensions}
\end{figure}
} 

\subsection{Common Priors}\label{sec:common-prior}
The common prior assumption simplifies the belief calculations in our hard
examples, but it does not play a critical role otherwise. In fact, we can argue
that similar to the directed links, imposition of common priors on the agents
simplifies their inference tasks. This is consistent with the fact that common
priors are crucial for reaching agreement \citep{aumann1976agreeing}.

We note that in the binary action model
the computations of agents with arbitrary priors can be reduced
to computations with uniform priors. One way to achieve this is as follows:
For each agent $j$ with a non-uniform prior $\nu_j$ we introduce
two auxiliary agents $\kappa_j$ and $\kappa'_j$ with uniform priors.
Agent $\kappa_j$ is observed by everyone, \emph{including} agent $j$,
while agent  $\kappa'_j$ is observed by everyone \emph{except} agent $j$ (see Figure \ref{fig:auxillary_non_uniform_priors}).
We then set the signal structures of agents $\kappa_j$ and $\kappa'_j$ such that
(cf.~\eqref{eq:11})
$\boldsymbol{\lambda}_{\kappa_j}(1) = \gamma_j = -\boldsymbol{\lambda}_{\kappa'_j}(0)$
and specify private signals $\mathbf{s}_{\kappa_j}=1$ and
$\mathbf{s}_{\kappa'_j}=0$. One can verify that:\endnote{In the context of Subsection~\ref{sec:average}, one might argue that in the absence of a common prior there is no fixed distribution of signals over which to obtain an average-case hardness result. Notwithstanding, the worst-case issue remains relevant because the observation history is now exponentially unlikely according to each agent's own prior.}
\begin{enumerate}[label=(\roman*)]
\item The signal of agent $\kappa_j$ effectively shifts the prior of agent $j$
  to $\nu_j$.
\item Since everyone observed $\kappa_j$,
the fact that the prior of agent $j$ has been shifted becomes
common knowledge.
\item No agent other than
  $j$ shifts their belief after observing both $\kappa_j$ and $\kappa'_j$.
\end{enumerate}

\subsection{I.I.D.~Signals}\label{sec:asymmtericlikelihoods}

Assuming that the private signals are (conditionally) i.i.d.~is common in social learning literature. It often simplifies the analysis and provides a useful approximation to study homogeneous populations. The signals in our reductions are not i.i.d., but this is only for convenience. In Appendix~\ref{sec:app:infomrativesignals}, we explain how to modify our proofs to work with i.i.d. signals.

To give a general idea, in each reduction there are two issues to deal with: First, the auxiliary agent $\kappa$ receives a special private signal with a distribution that is different form any other agents. In VERTEX-COVER, agent $\kappa$ receives a very strong private signal that induces a $\log$-belief ratio shift equivalent to $n_i = cn$ zero vertex agent signals for some constant $c>0$. Therefore, it is not surprising that we can replace $\kappa$ by $n_i$ agents with signal structure of vertex agents, all reporting zero private signals (cf.~Figure~\ref{fig:auxillary_n_i}). In EXACT-COVER, the auxiliary agent $\kappa$ receives a special signal that is twice as strong compared to the subset agents (its $\log$-likelihood ratio is twice the subset agent signals). We can use two i.i.d. signals to have the same effect, except that we need a mechanism to ensure that their signals agree (they are both zero, or both one). We can achieve this using five auxiliary agents as shown in Figure~\ref{fig:auxillary_iid_rx3}.

The second issue is that agent $i$, as well as the edge agents in VERTEX-COVER and the element agents in EXACT-COVER, do not receive private signals. This can be remedied by a similar idea as presented in Subsection~\ref{sec:common-prior}. In particular, we allow the agents to receive private signals which are then countervailed by matching opposite signals coming from auxiliary agents.

\subsection{Noisy Actions}\label{sec:noisy}

As discussed in Subsection~\ref{sec:literature},
\cite{aaronson2005complexity} shows that if two agents
decide to add noise to their exchanged opinions,
their rational beliefs can be approximated efficiently.
This is an interesting model in its own right: A typical approach to bounded
rationality needs to choose a rule for updating beliefs, and any such choice
is, to an extent, arbitrary. If, instead, it could be shown that ``noisy''
Bayesian updates are efficient, it would provide for an interesting alternative.

Notwithstanding, we show that adding noise does not change our hardness
results for network models. For concreteness, we focus on a particular modification of the
binary action model. However, we believe our ideas should work with most
other natural variants. More precisely, we consider the binary action model with an additional
parameter $0<\delta<1/2$. All the rules are the same except that
every time an agent broadcasts her opinion to the world, a glitch (bit flip) occurs
with probability $\delta$.

In other words, every time agent $i$ computes an action $\mathbf{a}_{i,t}=\mathbbm{1}(\boldsymbol{\mu}_{i,t}>1/2)$, its announced value ($\mathbf{a}'_{i,t}$) is flipped to $1-\mathbf{a}_{i,t}$, independently with probability $\delta$. We assume that all neighbors of $i$ observe the same action (as opposed to
flipping with probability $\delta$ independently for each neighbor).
Since the networks that we consider are acyclical, it does not matter if the
agents observe their own actions, i.e., if they learn that their actions were
flipped. As before, all these rules are common knowledge and the agents
estimate their beliefs ($\boldsymbol{\mu}_{i,t}$) using the Bayes rule.

In Appendix~\ref{app:noisy} we show that estimating beliefs in this model
is still NP-hard. The main idea
is that an agent in the noiseless binary action model
can be replaced with multiple copies
of noisy agents broadcasting the same action in such a way that the probability
of the transmission error is negligible compared to the other probabilities
that determine the computed beliefs (see Figure \ref{fig:auxillary_noisy}).

\section{Algorithms for Bayesian Choice}\label{sec:algorithms}

Refinement of information partitions with increasing observations is a key feature of rational learning problems and it is fundamental to major classical results that establish agreement (\cite{geanakoplos1982we}) or learning (\cite{BlackwellDubins62,LehrerSmorodinsky96}) among rational agents.
 
In the group decision setting, the list of possible signal profiles is regarded as the information set representing the current understanding of the agent about her environment, and the way additional observations are informative is by trimming the current information set and reducing the ambiguity in the set of initial signals that have caused the agent's history of past observations. Thereby, one can conceive a natural  method of computing agents' actions based on \emph{elimination of impossible signals}. By successively eliminating signals that are inconsistent with the new observations, we refine the partitions of the space of private signals, and at the same time, we keep track of the current information set that is consistent with the observations. As such, we refer to this approach as ``Elimination of Impossible Signals'' or EIS. We begin by presenting a recursive version (REIS), and study its iterative implementations (IEIS) afterwards.

To proceed, let $\overline{s} = (s_1,\ldots,s_n)\in\mathcal{S}_1\times\ldots\times\mathcal{S}_n$ be any profile of initial signals, and denote the set of all private signal profiles that agent $i$ regards as possible at time $t$, i.e. her information set at time $t$, by  $\bm{\mathcal{I}}_{i,t}\subset \mathcal{S}_1\times\ldots\times\mathcal{S}_n$; this random set is a function of the observed history $\mathbf{h}_{i,t}$ and is fully determined by the random profile of all private signals $\overline{\mathbf{s}} := (\mathbf{s}_1,\ldots,\mathbf{s}_n)$. Recall that the observation history $\mathbf{h}_{i,t}$ is defined as
$\{\mathbf{s}_i\} \cup \{\mathbf{a}_{j,\tau}$ for all $j\in\mathcal{N}_i$, and $\tau < t\}$. Hence, $\bm{\mathcal{I}}_{i,t}$ takes into account the neighboring actions at all times strictly less than $t$.

Starting from $\bm{\mathcal{I}}_{i,0} = \{\mathbf{s}_{i}\}\times\prod_{j\neq i}\mathcal{S}_j$, at every step $t>0$ agent $i$ removes those signal profiles in $\bm{\mathcal{I}}_{i,t-1}$ that are inconsistent with her observation history $\mathbf{h}_{i,t}$, and constructs a censured set of signal profiles $\bm{\mathcal{I}}_{i,t} \subset \bm{\mathcal{I}}_{i,t-1}$. Recall that ${\mathbb{P}}_{\theta}(\mathord{\cdot})$ is the joint distribution of the private signals of all agents. For each $i$ and $t$, the set of possible signals ($\bm{\mathcal{I}}_{i,t}$) is mapped to a Bayesian posterior (${\boldsymbol\mu}_{i,t}$) as follows:
\begin{align}
{\boldsymbol\mu}_{i,t}(\theta) = \frac{\sum_{\overline{s}\in\bm{\mathcal{I}}_{i,t}} {\mathbb{P}}_{\theta}(\overline{s})\nu(\theta)}{\sum_{{\theta'}\in\Theta}\sum_{\overline{s}\in\bm{\mathcal{I}}_{i,t}} {\mathbb{P}}_{\theta'}(\overline{s})\nu({\theta'})}. \label{eq:BayesianPosterior}
\end{align} The posterior belief, in turn, enables the agent to choose an optimal (myopic) action given her observations:
\begin{align}
 \mathbf{a}_{i,t} =  \argmax\limits_{a_i\in\mathcal{A}_i}\sum_{{\theta'}\in\Theta}u_i(a_i,{\theta'}){\boldsymbol\mu_{i,t}}({\theta'}). \label{eq:OptimalRecommendation}
\end{align}

It is convenient to define a function $\mathcal{A}_i$ that given a set of possible signal profiles $\mathcal{I} \subset \prod_{j=1}^{n}\mathcal{S}_j$ outputs the optimal action of agent $i$ as follows: 
\begin{align}
\mathcal{A}_i(\mathcal{I}) = \displaystyle \argmax_{a\in\mathcal{A}_i}\sum_{{\theta'}\in\Theta}u_i(a,{\theta'})\frac{\sum_{\overline{s}'\in\mathcal{I}} {\mathbb{P}}_{\theta'}(\overline{s}')\nu({\theta'})}{\sum_{{\theta''} \in \Theta}\sum_{\overline{s}'\in\mathcal{I}} {\mathbb{P}}_{\theta''}(\overline{s}')\nu({\theta''})}.
\label{eq:ActionsGivenSignalProfiles}
\end{align} Crucially, in addition to her own possible set $\boldsymbol{\mathcal{I}}_{i,t}$,
 agent $i$ keeps track of other agents' possible sets as well. Therefore, it is useful to consider the function ${\mathcal{I}}(j,t,\overline{s})$ that outputs the set of signal profiles that agent $j$ considers possible at time $t$ if the initial private signals are $\overline{s}$. Subsequently, the action that agent $j$ takes if the initial private signals are $\overline{s}$ is given by $\mathcal{A}_j(\mathcal{I}(j,t,\overline{s}))$. 

 Note that in the above notation, $\cup_{\overline{s}\in\boldsymbol{\mathcal{I}}_{i,t}} {\mathcal{I}}(j,t,\overline{s})$ is the set of all signal profiles that agent $i$ cannot yet conclude are rejected by agent $j$. Similarly, $\cup_{\overline{s}\in\boldsymbol{\mathcal{I}}_{i,t}}\mathcal{A}_{j}({\mathcal{I}}(j,t,\overline{s}))$ is the list of all possible actions that agent $j$ may currently take, from the viewpoint of agent $i$ (consistent with agent $i$'s observations so far). Given $\mathcal{A}_j({\mathcal{I}}(j,t,\overline{s}))$ for all $\overline{s} \in \bm{\mathcal{I}}_{i,t-1}$ and every $j \in \mathcal{N}_i$, agent $i$ can reject any $\overline{s}$ for which the observed neighboring action $\mathbf{a}_{j,t}$ does not agree with the simulated action: Reject any $\overline{s}$ such that $\mathbf{a}_{j,t} \neq \mathcal{A}_j(\mathcal{I}(j,t,\overline{s}))$  for some $j\in\mathcal{N}_i$.

The function ${\mathcal{I}}(i,t,\overline{s})$ can be defined recursively by listing all signal profiles $\overline{s}'$ that are consistent with $\overline{s}$, producing the same observations for agent $i$ up until time $t$. To check such consistencies one needs to make additional function calls of the form ${\mathcal{I}}(j,\tau,\overline{s}')$ for $j \in \mathcal{N}_i$ and $\tau< t$. We formalize this idea in Algorithm 0 by offering a recursive implementation for the elimination of impossible signals to compute Bayesian actions
(cf.~Table~\ref{varlist} for a summary of the notation).\endnote{We note that Algorithm~0 can be implemented to use space that is polynomial in the number of agents and time $t$ (assuming fixed state set $\Theta$, signal sets $\mathcal{S}_i$ and action sets $\mathcal{A}_i$). In the binary action model this matches our PSPACE-hardness results obtained in the follow-up paper~\cite{BayesPSPACE18}.}


\begin{table}
\caption{Notation for Bayesian group decision computations (Elimination of Impossible Signals)\label{varlist}}

{\begin{tabular}{l ||  p{0.7\textwidth}}
\hline
\hline
 \makecell{ \vspace{-10pt}\\$\overline{s} =$ \\ $ (s_1,s_2,\ldots,s_n)$} & \makecell[l]{ \vspace{-13pt}\\a profile of initial private signals.} \\
\hline
 \makecell{ \vspace{-2pt}\\${\bm{\mathcal{I}}}_{i,t}$} & the set of  all signal profiles that are deemed possible by agent $i$, given her observations up until time $t$. \\
\hline
 \makecell{ \vspace{-2pt}\\${\mathcal{I}}(j,t,\overline{s})$} &  the set of all signal profiles that are deemed possible by agent $j$ at time $t$, if the initial signals of all agents are prescribed according to $\overline{s}$. \\
\hline
\makecell{ \vspace{-2pt}\\ $A_{j}({\mathcal{I}}(j,t,\overline{s}))$} & the computed action of agent $j$ at time $t$, if  the initial signals of all agents are prescribed according to $\overline{s}$. \\
\hline
\hline
\end{tabular}}
\end{table}


\noindent
 \begingroup\fboxsep=10pt
 \begin{center}
\fbox{\begin{minipage}{0.9\textwidth}
\parbox{\textwidth}{{\bf Algorithm 0: RECURSIVE-EIS $(i,t)$}\vspace{0.0cm}\\ 
    {\bf Input:} Graph $\mathcal{G}$, set of possible signal profiles $\bm{\mathcal{I}}_{i,t}$, and neighboring actions $\mathbf{a}_{j,t}, j \in\mathcal{N}_i$ \vspace{0.0cm}\\
    {\bf Output:} Bayesian action $\mathbf{a}_{i,t+1}$\vspace{0.05cm}
    \begin{enumerate}
	    \item Initialize $\bm{\mathcal{I}}_{i,t+1} = \bm{\mathcal{I}}_{i,t}$.
	    \item For all  $\overline{s}\in\bm{\mathcal{I}}_{i,t+1}$, do:
	    \begin{itemize}
	        \item For all $j\in\mathcal{N}_i$, if $\mathbf{a}_{j,t} \neq \mathcal{A}_{j}(\mathcal{I}(j,t,\overline{s}))$, then set $ \bm{\mathcal{I}}_{i,t+1} = \bm{\mathcal{I}}_{i,t+1}\setminus \{ \overline{s} \}$.
	    \end{itemize}
	    \item  $\mathbf{a}_{i,t+1} = \mathcal{A}_i(\bm{\mathcal{I}}_{i,t+1})$.
	\end{enumerate}

{\bf Function} $\mathcal{I}(i,t,\overline{s}):$
\begin{itemize}
        \item If $t = 0$, then set $I = \{{s}_{i}\}\times\prod_{j\neq i}\mathcal{S}_j$
        \item else if $t > 0$:
        \begin{enumerate}
               \item Initialize $I = \varnothing$.
               \item For all $\overline{s}'\in \mathcal{S}_1\times\ldots\times\mathcal{S}_n$, do:
               \begin{itemize}
                       \item If \emph{Consistent}$(i,t,\overline{s},\overline{s}')$, then set $I = I \cup \{\overline{s}'\}$.
                \end{itemize}
        \end{enumerate}
\end{itemize}        
{\bf return} $I$

{\bf Function} \emph{Consistent}$(i,t,\overline{s},\overline{s}')$:
\begin{enumerate}
    \item Initialize \texttt{is\_consistent} $=$ \texttt{True}.
    \item For all $\tau < t$ and $j \in \mathcal{N}_i$, do: 
    \begin{itemize}
        \item If $\mathcal{A}_j(\mathcal{I}(j,\tau,\overline{s})) \neq \mathcal{A}_j(\mathcal{I}(j,\tau,\overline{s}'))$, then \texttt{is\_consistent} $=$ \texttt{False}.
    \end{itemize}
\end{enumerate}

{\bf return} \texttt{is\_consistent}

}
\end{minipage}}
\end{center}
\endgroup
\vspace{10pt}

In Subsection \ref{sec:generalStrcutures}, we describe an iterative implementation of elimination of impossible signals (IEIS). The IEIS calculations scale exponentially with the network size; this is true, in general, with the exception of some  densely connected networks where agents have direct access to all the observations of their neighbors. We expand on this special case (called transitive networks) in Subsection \ref{sec:POSETs}. Finally, in Subsection~\ref{sec:algobelief}
we discuss the revealed beliefs case and identify additional network structures for which Bayesian calculations simplify, allowing for efficient Bayesian belief exchange.

\subsection{Iterative Elimination of Impossible Signals (IEIS)}\label{sec:generalStrcutures}
 
To proceed, we denote $\mathcal{N}^{\tau}_i$ as the $\tau$-th order neighborhood of agent $i$ comprising entirely of those agents who are at distance $\tau$ from agent $i$; in particular, $\mathcal{N}^{1}_{i} = \mathcal{N}_{i}$, and we use the convention $\mathcal{N}^{0}_{i} = \{{i}\}$. We further denote $\bar{\mathcal{N}}_i^{t}:=\cup_{\tau=0}^{t}\mathcal{N}^{\tau}_i$ as the set of all agents who are within distance $t$ of or closer to agent $i$; we sometimes refer to $\bar{\mathcal{N}}_i^{t}$ as her ego-net of radius $t$.

At time zero, agent $i$ initializes her list of possible signals $\bm{\mathcal{I}}_{i,0} = \{\mathbf{s}_{i}\}\times\prod_{j\neq i}\mathcal{S}_j$. At time $t$, she has access to $\bm{\mathcal{I}}_{i,t}$, the list of possible signal profiles that are consistent with  her observations so far, as well as all signal profiles that she thinks each of the other agents would regard as possible conditioned on any profile of initial signals: $\mathcal{I}(j,t-\tau,\overline{s})$ for $\overline{s}\in \mathcal{S}_1\times\ldots\times\mathcal{S}_n$, $j\in \mathcal{N}^{\tau}_i$, and $\tau \in [t]:=\{1,2,\ldots,t\}$. Given the newly obtained information, which constitute her observations of the most recent neighboring actions $\mathbf{a}_{j,t}$, $j \in \mathcal{N}_{i}$, she refines $\bm{\mathcal{I}}_{i,t}$ to $\bm{\mathcal{I}}_{i,t+1}$ and updates her belief and actions accordingly, cf. \eqref{eq:BayesianPosterior} and \eqref{eq:OptimalRecommendation}. This is achieved as follows (we use $\mbox{dist}(j,i)$ to denote the length of the shortest path connecting $j$ to $i$):

\vspace{10pt}

\noindent
 \begingroup\fboxsep=10pt
 \begin{center}
\fbox{
\begin{minipage}{0.8\textwidth}
\parbox{\textwidth}{{\bf Algorithm 1: IEIS $(i,t)$}\vspace{0.0cm}\\ 
    {\bf Input:} Graph $\mathcal{G}$, \\ set of possible signal profiles $\bm{\mathcal{I}}_{i,t}$, $\mathcal{I}(j,\tau,\overline{s})$, for all $\overline{s}$, $\tau \in [t-\mbox{dist}(j,i)]$, $j\in \bar{\mathcal{N}}^{t}_i$,\\ and neighboring actions $\mathbf{a}_{j,t}, j \in\mathcal{N}_i$ \vspace{0.0cm}\\
    {\bf Output:} Bayesian action $\mathbf{a}_{i,t+1}$\vspace{0.05cm}
 
\begin{itemize}[leftmargin=*]
\item SIMULATE:\\ 
For all  $\overline{s}:=(s_1 , \ldots , s_n) \in \mathcal{S}_1\times\ldots\times\mathcal{S}_n$, do: 
	\begin{enumerate}
    \item For $j\in \mathcal{N}^{t+1}_i$, initialize  $\mathcal{I}(j,0,\overline{s}) = \{s_j\}\times\prod_{k\neq j}\mathcal{S}_k$.
	\item For $\tau=t, t-1, \ldots, 1$, do: 
	\begin{enumerate}
	    \item For $j\in \mathcal{N}^{\tau}_i$, do:
	    \begin{enumerate}
	    \item Initialize $\mathcal{I}(j,t+1-\tau,\overline{s}) = \mathcal{I}(j,t-\tau,\overline{s})$.
	    \item For $\overline{s}' \in \mathcal{I}(j,t+1-\tau,\overline{s})$ do:
		\begin{itemize}
		\item  For all $k\in\mathcal{N}_j$, if $\mathcal{A}_k(\mathcal{I}(k,t-\tau,\overline{s}')) \neq \mathcal{A}_k(\mathcal{I}(k,t-\tau,\overline{s}))$, \\ then set $ \mathcal{I}(j,t+1-\tau,\overline{s}) = \mathcal{I}(j,t+1-\tau,\overline{s})\setminus \{ \overline{s}' \}$.
		\end{itemize}
	\end{enumerate}
	\end{enumerate}
	\end{enumerate}	
\item UPDATE: 
	\begin{enumerate}[leftmargin=*]
	    \item Initialize $\bm{\mathcal{I}}_{i,t+1} = \bm{\mathcal{I}}_{i,t}$.
	    \item For all  $\overline{s}\in\bm{\mathcal{I}}_{i,t+1}$, do:
	    \begin{itemize}
	        \item For all $j\in\mathcal{N}_i$, if $\mathbf{a}_{j,t} \neq \mathcal{A}_j(\mathcal{I}(j,t,\overline{s}))$, then set $ \bm{\mathcal{I}}_{i,t+1} = \bm{\mathcal{I}}_{i,t+1}\setminus \{ \overline{s} \}$.
	    \end{itemize}
	    \item Set $\mathbf{a}_{i,t+1} = \mathcal{A}_i(\bm{\mathcal{I}}_{i,t+1})$.
	\end{enumerate}
\end{itemize}
}
\end{minipage}}
\end{center}
\endgroup
\vspace{10pt} 

\vspace{10pt}

Note that in the ``SIMULATE'' part of the IEIS Algorithm, we make no use of the observations of agent $i$. This step amounts to simulating the network at all signal profiles. It is implemented such that the computations at time $t$ are based on what was computed for making decisions prior to time $t$. In the ``UPDATE'' part, we compare the most recently observed actions of neighbors with their simulated actions for each signal profile in $\bm{\mathcal{I}}_{i,t}$ to detect and eliminate the impossible ones. To evaluate the possibility of a signal profile using IEIS, agent $i$ may need to consider actions that other agents could have taken in signal profiles that she has already rejected. In particular, simulating the network at all possible profiles of agent $i$ at time $t$, i.e. at all $ \overline{s}\in\bm{\mathcal{I}}_{i,t}$, is not enough to evaluate the condition, $\mathcal{A}_k(\mathcal{I}(k,t-\tau,\overline{s}')) \neq \mathcal{A}_k(\mathcal{I}(k,t-\tau,\overline{s}))$, at step 2(a)ii of Algorithm~1--SIMULATE, since $\overline{s}'$ may not be included in $\bm{\mathcal{I}}_{i,t}$.

In Appendix \ref{App:BAYESGROUPCOMP} we describe the complexity of the computations that the agent should undertake using IEIS at any time $t$ in order to calculate her posterior probability ${\boldsymbol\mu}_{i,t+1}$ and Bayesian decision $\mathbf{a}_{i,t+1}$ given all her observations up to time $t$. Subsequently, we prove that:

\begin{theorem}[Complexity of IEIS]  Consider a network of size $n$ with $m$ states, and let $M$ and $A$ denote the maximum cardinality of the signal and action spaces $(m:=$\emph{card}$(\Theta)$, $M = \max_{k\in[n]}$\emph{card}$(\mathcal{S}_k)$, and $A = \max_{k\in[n]}$\emph{card}$(\mathcal{A}_k))$. The IEIS algorithm has $O(n^2 M^{2n-1} m A)$ running time, which given the private signal of agent $i$ and the previous actions of her neighbors $\{\mathbf{a}_{j,\tau}:j\in\mathcal{N}_i,\tau< t\}$ in any network structure, outputs $\mathbf{a}_{i,t}$, the Bayesian action of agent $i$ at a fixed time $t$. 
 \end{theorem}

\subsection{IEIS over Transitive Structures}\label{sec:POSETs}

We now shift focus to the special case of transitive networks, defined below. 
 
\begin{definition}[Transitive Networks]\label{def:POSETnets} We call a network structure transitive if the directed neighborhood relationship between its nodes satisfies the reflexive and transitive properties. In particular, the transitive property implies that anyone whose actions indirectly influence the observations of agent $i$ is also directly observed by her, i.e. any neighbor of a neighbor of agent $i$  is a neighbor of agent $i$ as well. \end{definition}

In such structures, any agent whose actions indirectly influence the observations of agent $i$ is also directly observed by her. This special structure of transitive networks mitigates the issue of hidden observations, and as a result, Bayesian inference in a transitive structure is significantly less complex.

After initializing $\bm{\mathcal{S}}_{j,0} = \mathcal{S}_j$ and $\bm{\mathcal{I}}_{i,0} = \{\mathbf{s}_{i}\}\times\prod_{j\in\mathcal{N}_ {i}}{\bm{\mathcal{S}}}_{j,0}$, agent $i$ needs only to keep track of ${\bm{\mathcal{S}}}_{j,t} \subseteq \mathcal{S}_j$ for all $j\in\mathcal{N}_i$ (cf. Table \ref{varlist2}). This is because, in transitive structures, the list of possible signal profiles decomposes: $\bm{\mathcal{I}}_{i,t} = \{\mathbf{s}_i\}\times\prod_{j\in\mathcal{N}_i}\bm{\mathcal{S}}_{j,t}$. Updating in transitive structures is achieved by incorporating $\mathbf{a}_{j,t}$ for each $j\in\mathcal{N}_i$ individually, and transforming the respective $\bm{\mathcal{S}}_{j,t}$ into $\bm{\mathcal{S}}_{j,t+1}$. This updating procedure is formalized in Algorithm~2.

\begin{table}
\caption{Notation for Computations in Transitive Networks \label{varlist2}}
{\begin{tabular}{l ||  p{0.75\textwidth}}
\hline
\hline

\makecell{ \vspace{-5pt}\\ ${\bm{\mathcal{S}}}_{i,t}$} &  \makecell[l]{ \vspace{-10pt}\\the list of  all private signals that are deemed possible for agent $i$ at time $t$,\\ by an agent who has observed her actions in a transitive network structure\\ up until time $t$.} \\
\hline
\makecell{ \vspace{-10pt}\\ $\bm{\mathcal{I}}_{i,t}(s_i) = \{s_i\}\times$} \\ $ \prod_{j\in\mathcal{N}_i}{\bm{\mathcal{S}}_{j,t}}$ &  \makecell[l]{ \vspace{-20pt}\\   the list of neighboring signal profiles that are deemed possible by agent $i$,\\ given her observations of their actions up until time $t$ conditioned on own\\ private signal being $s_i$. }\\
\hline
\hline
\end{tabular}}
\end{table}

\vspace{10pt}

\noindent
 \begingroup\fboxsep=10pt
 \begin{center}
\fbox{
\begin{minipage}{0.7\textwidth}
\parbox{\textwidth}{{\bf Algorithm 2: IEIS-TRANSITIVE $(i,t)$}\vspace{0.0cm}\\ 
    {\bf Input:} Transitive graph $\mathcal{G}$, set of possible signal profiles $\bm{\mathcal{S}}_{j,t}$, $\forall j\in {\mathcal{N}}_i$, and neighboring actions $\mathbf{a}_{j,t}, j \in\mathcal{N}_i$. \vspace{0.0cm}\\
    {\bf Output:} Bayesian action $\mathbf{a}_{i,t+1}$\vspace{0.05cm}
\begin{enumerate}[leftmargin=*]
\item For all $j\in \mathcal{N}_i$, do:
	\begin{enumerate}
	\item Initialize $\bm{\mathcal{S}}_{j,t+1} = \bm{\mathcal{S}}_{j,t}$. 
	\item For all  $s_j \in \bm{\mathcal{S}}_{j,t+1}$, do: 
	\begin{enumerate}
	\item Set $\bm{\mathcal{I}}_{j,t}(s_j) = \{s_j\}\times\prod_{k\in\mathcal{N}_j}{\bm{\mathcal{S}}_{k,t}}$.
	\item If $\mathbf{a}_{j,t}\neq \mathcal{A}_j(\bm{\mathcal{I}}_{j,t}(s_j))$,  then set $\bm{\mathcal{S}}_{j,t+1} = \bm{\mathcal{S}}_{j,t+1}\setminus \{ {s}_j \}$.
	\end{enumerate}
	\end{enumerate}
\item Update $\bm{\mathcal{I}}_{i,t+1} = \{\mathbf{s}_i\}\times\prod_{j\in\mathcal{N}_i}{\bm{\mathcal{S}}_{j,t+1}}$.
\item Set $\mathbf{a}_{i,t+1} = \mathcal{A}_i(\bm{\mathcal{I}}_{i,t+1})$.
\end{enumerate}
}
\end{minipage}}
\end{center}
\endgroup
\vspace{10pt} 

In Appendix \ref{App:complexityOFposets}, we determine the computational complexity of the IEIS-TRANSITIVE algorithm as follows:

\begin{theorem} [Efficient Bayesian group decisions in transitive structures] Consider a network of size $n$ with $m$ states, and let $M$ and $A$ denote the maximum cardinality of the signal and action spaces $(m:=$\emph{card}$(\Theta)$, $M = \max_{k\in[n]}$\emph{card}$(\mathcal{S}_k)$, and $A = \max_{k\in[n]}$\emph{card}$(\mathcal{A}_k))$. There exists an algorithm with running time $O(A m n^2 M^2)$ which given the private signal of agent $i$ and the previous actions of her neighbors $\{\mathbf{a}_{j,\tau}:j\in\mathcal{N}_i,\tau< t\}$ in any transitive network, outputs $\mathbf{a}_{i,t}$, the Bayesian action of agent $i$ at time $t$. 
\end{theorem} 

\subsection{Algorithms for Beliefs}
\label{sec:algobelief}
  
In general, GROUP-DECISION with revealed beliefs is a hard problem per Theorem \ref{THM:NPHARDNESSBELIEFS}. Here, we introduce a structural property of the networks, called ``transparency'', which leads to efficient belief calculations in the revealed belief model. Recall that the $t$-radius ego-net of agent $i$, $\bar{\mathcal{N}}_i^{t}$, is the set of all agents who are within distance $t$ of or closer to agent $i$. In a transparent network, the belief of every agent at time $t$ aggregates the likelihoods of all private signals in their $t$-radius ego-net:

\begin{definition}[Transparency]\label{def:transparency} The graph structure $\mathcal{G}$ is transparent if for all agents $i \in [n]$ and all times $t$ we have that: $\boldsymbol{\phi}_{i,t} = \sum_{j\in\bar{\mathcal{N}}_i^{t}}\boldsymbol{\lambda}_j$,  for any choice of signal structures and all possible initial signals. Moreover, we call $\mathcal{G}$ transparent to agent $i$ at time $t$, if for all $j\in\mathcal{N}_i$ and every $\tau\leq t-1$ we have that: $\boldsymbol{\phi}_{j,\tau} = \sum_{k\in\bar{\mathcal{N}}_j^{\tau}}\boldsymbol{\lambda}_k$,  for any choice of signal structures and all possible initial signals.
\end{definition}

In any graph structure, the initial belief exchange between the agents reveals the likelihoods of the private signals in the neighboring agents (see Example \ref{example:initial_rounds} and equation \eqref{eq:beliefexchange} therein). Hence, from her observations of the beliefs of her neighbors at time zero, agent $i$ learns all that she needs to know regarding their private signals: 

\begin{corollary}[Transparency at time one] All graphs are transparent at time one.
\end{corollary}

However, the future neighboring beliefs (at time two and beyond) are ``less \emph{transparent}'' when it comes to reflecting the neighbors' knowledge of other private signals that are received throughout the network. In particular, the time one beliefs of the neighbors $\boldsymbol{\phi}_{j,1}, j\in\mathcal{N}_i$ are given by  $\boldsymbol{\phi}_{j,1} = \sum_{k\in\bar{\mathcal{N}}^{1}_j}\boldsymbol{\lambda}_k$; hence, from observing the time one belief of a neighbor, agent $i$ would only get to know $\sum_{k\in\mathcal{N}_j}\boldsymbol{\lambda}_k$, rather than the individual values of $\boldsymbol{\lambda}_k$ for each $k\in\mathcal{N}_j$.\endnote{This is a fundamental aspect of inference problems in observational learning (in learning from other actors): similar to responsiveness that \cite{ali2018herding} defines as a property of the utility functions to determine whether players' beliefs can be inferred from their actions, \emph{transparency} in our belief exchange setup is defined as a property of the graph structure (see Remark \ref{rem:transparency} on why transparency is a structural property) which determines to what extent other players' private signals can be inferred from observing the neighboring beliefs.}

\begin{remark}[Transparency, statistical efficiency, and impartial inference]\label{rem:transparency} Such agents $j$ whose beliefs satisfy the equation in Definition \ref{def:transparency} at some time $\tau$ are said to hold a \emph{transparent} or \emph{efficient} belief; the latter signifies the fact that such a belief coincides with the Bayesian posterior if agent $j$ were given direct access to the private signals of every agent in $\bar{\mathcal{N}}_j^{\tau}$. This is indeed the best possible (or statistically efficient) belief that agent $j$ can hope to form given the information available to her at time $\tau$. The same connection to the statistically efficient beliefs arise in the work of  \cite{eyster2014extensive} who formulate the closely related concept of ``impartial inference'' in a model of sequential decisions by different players in successive rounds; accordingly, impartial inference ensures that the full informational content of all signals that influence a player's beliefs can be extracted and players can fully (rather than partially) infer their predecessors' signals. In other words, under impartial inference, players' immediate predecessors provide ``sufficient statistics'' for earlier movers that are indirectly observed \cite[Section 3]{eyster2014extensive}.  Last but not least, it is worth noting that statistical efficiency or impartial inference are properties of the posterior beliefs, and as such the signal structures may be designed so that  statistical efficiency or impartial inference hold true for a particular problem setting; on the other hand, transparency is a structural property of the network and would hold true for any choice of  signal structures and all possible initial signals.  \end{remark}

Our next example helps clarify the concept of transparency as a structural graph property, and its relation to Bayesian belief computations.   

\begin{example}[Transparent Structures]\label{example:trasnparency}

{\begin{figure}[ht]
\centering
\begin{subfigure}[b]{0.15\textwidth}
   \includegraphics[width=1\linewidth]{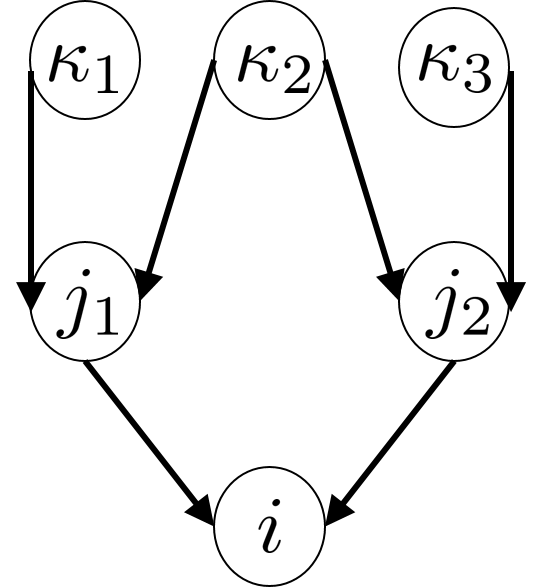}
   \caption{ }
   \label{fig:NonTransparenta}
\end{subfigure}~
\begin{subfigure}[b]{0.15\textwidth}
  \includegraphics[width=1\linewidth]{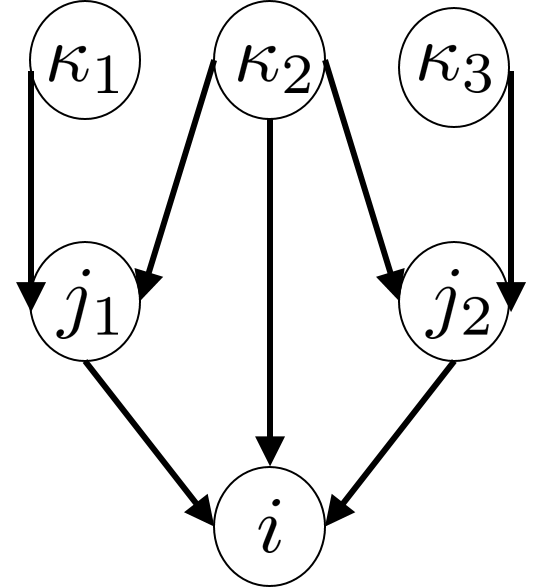}
\caption{ }
\label{fig:Transparent1a}
\end{subfigure}~
\begin{subfigure}[b]{0.15\textwidth}
  \includegraphics[width=1\linewidth]{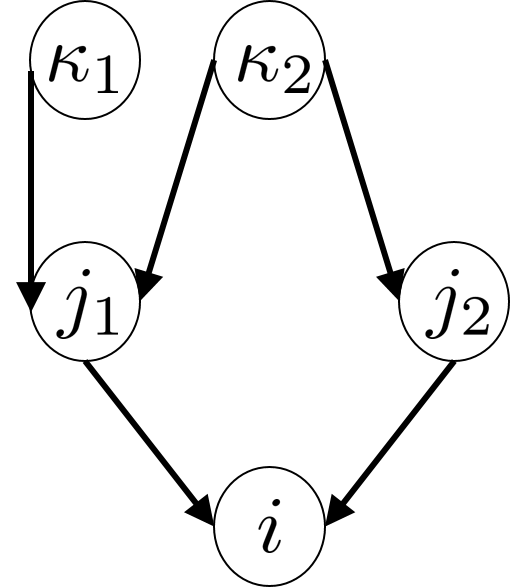}
\caption{ }
\label{fig:Transparent2a}
\end{subfigure}~
\begin{subfigure}[b]{0.25\textwidth}
  \includegraphics[width=1\linewidth]{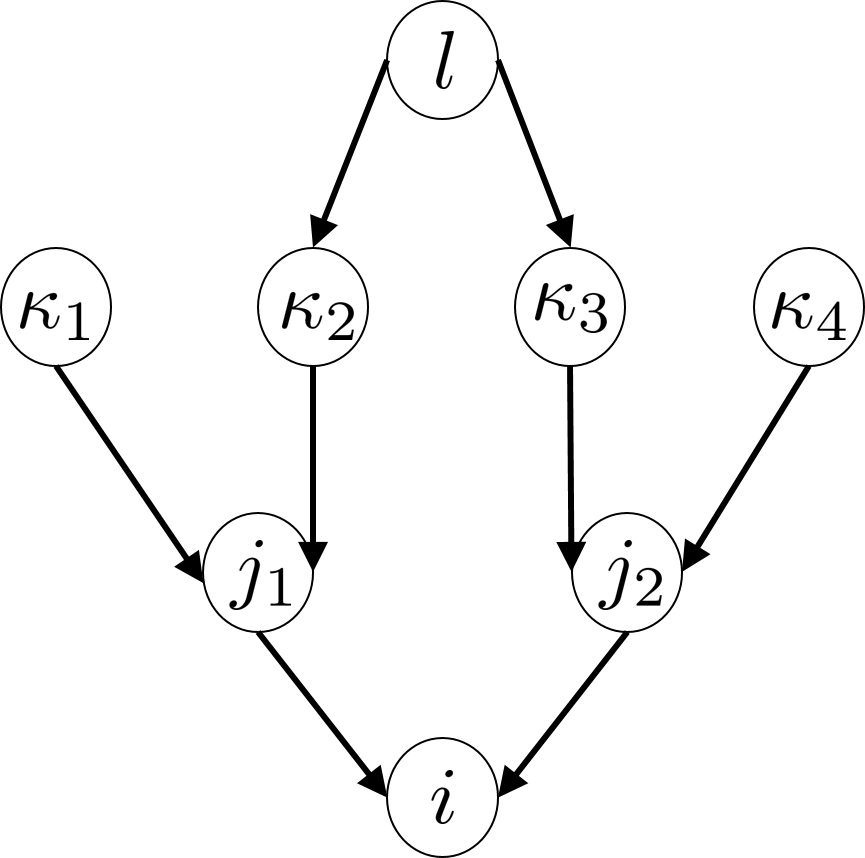}
\caption{ }
\label{fig:Transparent3a}
\end{subfigure}
\caption{Structures (\subref{fig:Transparent1a}-\subref{fig:Transparent3a}) are transparent, but (\subref{fig:NonTransparenta}) is not.}
\label{fig:transparency}
\end{figure}
}

Figure \ref{fig:transparency} illustrates cases of transparent and nontransparent structures. All structures except (\subref{fig:NonTransparenta}) are transparent. To see how the transparency is violated in (\subref{fig:NonTransparenta}), consider the beliefs of agent $i$: 
\begin{align}
{\boldsymbol{\phi}}_{i,0} &= {\boldsymbol{\lambda}}_{i},\\ 
{\boldsymbol{\phi}}_{i,1} &= {\boldsymbol{\lambda}}_{i} + {\boldsymbol{\lambda}}_{j_1} + {\boldsymbol{\lambda}}_{j_2}.    
\end{align} At time two, agent one observes the following reports: 
\begin{align}
{\boldsymbol{\phi}}_{j_1,1} &= {\boldsymbol{\lambda}}_{j_1} + {\boldsymbol{\lambda}}_{\kappa_1} + {\boldsymbol{\lambda}}_{\kappa_2},\\ {\boldsymbol{\phi}}_{j_2,1} &= {\boldsymbol{\lambda}}_{j_2} + {\boldsymbol{\lambda}}_{\kappa_2} + {\boldsymbol{\lambda}}_{\kappa_3}.    
\end{align} Knowing ${\boldsymbol{\phi}}_{j_1,0} = {\boldsymbol{\lambda}}_{j_1}$ and ${\boldsymbol{\phi}}_{j_2,0} = {\boldsymbol{\lambda}}_{j_2}$ she can infer the values of the two sub-sums  ${\boldsymbol{\lambda}}_{\kappa_1} + {\boldsymbol{\lambda}}_{\kappa_2}$ and ${\boldsymbol{\lambda}}_{\kappa_2} + {\boldsymbol{\lambda}}_{\kappa_3}$, but there is no way for her to infer their total sum $ {\boldsymbol{\lambda}}_{j_1} +  {\boldsymbol{\lambda}}_{j_2} + {\boldsymbol{\lambda}}_{\kappa_1} + {\boldsymbol{\lambda}}_{\kappa_2}+{\boldsymbol{\lambda}}_{\kappa_3}$. Agent $i$ cannot hold a belief that efficiently aggregates all private signals at time two; hence, the first structure is not transparent. Here, it is instructive to exactly characterize the non-transparent Bayesian posterior belief of agent $i$ at time two. At time two, agent $i$ can determine the sub-sum $\boldsymbol{\lambda}_i + \boldsymbol{\lambda}_{j_1} +  \boldsymbol{\lambda}_{j_2}$ and her belief would involve a search only over the profile of the signals of the remaining agents $(s_{\kappa_1},s_{\kappa_2},s_{\kappa_3})$. At time two, she finds all $(s_{\kappa_1},s_{\kappa_2},s_{\kappa_3})$ that agree with the additionally inferred sub-sums  $\boldsymbol{\lambda}_{\kappa_1} + \boldsymbol{\lambda}_{\kappa_2}$ and $\boldsymbol{\lambda}_{\kappa_2} +  \boldsymbol{\lambda}_{\kappa_3}$. If we use $ \bm{\mathcal{I}}_{i,2}$ to denote the set of all such triplets of feasible signals $(s_{\kappa_1},s_{\kappa_2},s_{\kappa_3})$, then  we can express $\boldsymbol\phi_{i,2}$ as follows:
\begin{align}
&\boldsymbol{\phi}_{i,2} = \boldsymbol{\lambda}_i + \boldsymbol{\lambda}_{j_1} +  \boldsymbol{\lambda}_{j_2} + \log\frac{\sum_{(s_{\kappa_1},s_{\kappa_2},s_{\kappa_3}) \in \bm{\mathcal{I}}_{i,2}}{\mathbb{P}}_{\kappa_1,\theta_2}(s_{\kappa_1}){\mathbb{P}}_{\kappa_2,\theta_2}(s_{\kappa_2}){\mathbb{P}}_{\kappa_3,\theta_2}(s_{\kappa_3})}{\sum_{(s_{\kappa_1},s_{\kappa_2},s_{\kappa_3}) \in \bm{\mathcal{I}}_{i,2}}{\mathbb{P}}_{\kappa_1,\theta_1}(s_{\kappa_1}){\mathbb{P}}_{\kappa_2,\theta_1}(s_{\kappa_2}){\mathbb{P}}_{\kappa_3,\theta_1}(s_{\kappa_3})}, \label{eq:highlynonlinerbeliefs}
\end{align}
 where  
 \begin{align} \bm{\mathcal{I}}_{i,2} = \{(s_{\kappa_1},s_{\kappa_2},s_{\kappa_3}): & \log\frac{{\mathbb{P}}_{\kappa_1,\theta_2}(s_{\kappa_1})}{{\mathbb{P}}_{\kappa_1,\theta_1}(s_{\kappa_1})}+  \log\frac{{\mathbb{P}}_{\kappa_2,\theta_2}(s_{\kappa_2})}{{\mathbb{P}}_{\kappa_2,\theta_1}(s_{\kappa_2})} = \boldsymbol{\lambda}_{\kappa_1}   +  \boldsymbol{\lambda}_{\kappa_2},  \mbox{ and }  \\ & \log\frac{{\mathbb{P}}_{\kappa_1,\theta_2}(s_{\kappa_1})}{{\mathbb{P}}_{\kappa_1,\theta_1}(s_{\kappa_1})} + \log\frac{{\mathbb{P}}_{\kappa_3,\theta_2}(s_{\kappa_3})}{{\mathbb{P}}_{\kappa_3,\theta_1}(s_{\kappa_3})} = \boldsymbol{\lambda}_{\kappa_2} +  \boldsymbol{\lambda}_{\kappa_3}\}.
 \end{align} 

We now move to the next structure (\subref{fig:Transparent1a}).  The ambiguity in determining ${\boldsymbol{\lambda}}_{\kappa_1} + {\boldsymbol{\lambda}}_{\kappa_2}+{\boldsymbol{\lambda}}_{\kappa_3}$ is resolved in (\subref{fig:Transparent1a}) by simply adding a direct link so that agent $\kappa_2$ is directly observed by agent $i$. Subsequently, agent $i$ holds an efficient posterior belief at time two: ${\boldsymbol{\phi}}_{i,2} =  {\boldsymbol{\lambda}}_{i} + {\boldsymbol{\lambda}}_{j_1} +  {\boldsymbol{\lambda}}_{j_2} + {\boldsymbol{\lambda}}_{\kappa_1} + {\boldsymbol{\lambda}}_{\kappa_2}+{\boldsymbol{\lambda}}_{\kappa_3}$.

In (\subref{fig:Transparent2a}), agent $i$ observes the following reports of her neighbors:
\begin{align}
    {\boldsymbol{\phi}}_{j_1,0} & = {\boldsymbol{\lambda}}_{j_1}, \\
    {\boldsymbol{\phi}}_{j_2,0} & = {\boldsymbol{\lambda}}_{j_2}, \\
    {\boldsymbol{\phi}}_{j_1,1} & =  {\boldsymbol{\lambda}}_{j_1} + {\boldsymbol{\lambda}}_{\kappa_1} + {\boldsymbol{\lambda}}_{\kappa_2},
\end{align} and can use these observations at time two, to solve for the sum of $\log$-likelihood ratios of private signals of everybody:
\begin{align}
    {\boldsymbol{\phi}}_{i,2} &= {\boldsymbol{\lambda}}_{i} + {\boldsymbol{\phi}}_{j_1,1}+ {\boldsymbol{\phi}}_{j_2,0} \\ 
    & = {\boldsymbol{\lambda}}_{i} + {\boldsymbol{\lambda}}_{j_1} +  {\boldsymbol{\lambda}}_{j_2} + {\boldsymbol{\lambda}}_{\kappa_1} + {\boldsymbol{\lambda}}_{\kappa_2} 
\end{align}

 Structure (\subref{fig:Transparent3a}) is also transparent. At time two, agent $i$ observes ${\boldsymbol{\phi}}_{j_1,1} = {\boldsymbol{\lambda}}_{j_1} + {\boldsymbol{\lambda}}_{\kappa_1} + {\boldsymbol{\lambda}}_{\kappa_2}$ and $ {\boldsymbol{\phi}}_{j_2,1} = {\boldsymbol{\lambda}}_{j_2} + {\boldsymbol{\lambda}}_{\kappa_3} + {\boldsymbol{\lambda}}_{\kappa_4}$, in addition to her own private signal ${\boldsymbol{\lambda}}_{i}$. Her belief at time two is given by: 
 \begin{align}
    {\boldsymbol{\phi}}_{i,2} &= {\boldsymbol{\lambda}}_{i} + {\boldsymbol{\phi}}_{j_1,1}+ {\boldsymbol{\phi}}_{j_2,1} \\ 
    & = {\boldsymbol{\lambda}}_{i} + {\boldsymbol{\lambda}}_{j_1} + {\boldsymbol{\lambda}}_{j_2} + {\boldsymbol{\lambda}}_{\kappa_1} + {\boldsymbol{\lambda}}_{\kappa_2} + {\boldsymbol{\lambda}}_{\kappa_3} + {\boldsymbol{\lambda}}_{\kappa_4}.
 \end{align} A time three, agent $i$ adds ${\boldsymbol{\phi}}_{j_1,2} = {\boldsymbol{\phi}}_{j_1,2} + {\boldsymbol{\lambda}}_{l} = {\boldsymbol{\lambda}}_{j_1} + {\boldsymbol{\lambda}}_{\kappa_1} + {\boldsymbol{\lambda}}_{\kappa_2} + {\boldsymbol{\lambda}}_{l}$ to her observations and her belief at time three is give by:
 \begin{align}
     {\boldsymbol{\phi}}_{i,3} &= {\boldsymbol{\lambda}}_{i} + {\boldsymbol{\phi}}_{j_1,1}+ {\boldsymbol{\phi}}_{j_2,1}+ ({\boldsymbol{\phi}}_{j_1,2} - {\boldsymbol{\phi}}_{j_1,1}) \\ & = {\boldsymbol{\lambda}}_{i} + {\boldsymbol{\lambda}}_{j_1} +  {\boldsymbol{\lambda}}_{j_2} + {\boldsymbol{\lambda}}_{\kappa_1} + {\boldsymbol{\lambda}}_{\kappa_2}+ {\boldsymbol{\lambda}}_{\kappa_3} + {\boldsymbol{\lambda}}_{\kappa_4} + {\boldsymbol{\lambda}}_{l}.
 \end{align} This example illustrates a case where an agent learns the sum of $\log$-likelihood ratios of signals of agents in her higher-order neighborhoods even though she cannot determine each  $\log$-likelihood ratio individually. In structure (\subref{fig:Transparent3a}), agent $i$ learns  $\{{\boldsymbol\lambda}_{i} ,{\boldsymbol\lambda}_{j_1} , {\boldsymbol\lambda}_{j_2} , {\boldsymbol\lambda}_{\kappa_1} + {\boldsymbol\lambda}_{\kappa_2} , {\boldsymbol\lambda}_{\kappa_3} + {\boldsymbol\lambda}_{\kappa_4} ,  {\boldsymbol\lambda}_{l}\}$, and in particular, she can determine the total sum of $\log$-likelihood ratios of all of the signals in her extended neighborhood, but she never learns the values of the individual $\log$-likelihood ratios $\{{\boldsymbol\lambda}_{\kappa_1} , {\boldsymbol\lambda}_{\kappa_2} , {\boldsymbol\lambda}_{\kappa_3} , {\boldsymbol\lambda}_{\kappa_4}\}$. \QEDB

\end{example}

The following is a sufficient graphical condition for agent $i$ to hold an efficient (transparent) belief at time $t$: there are no agents $k \in \bar{\mathcal{N}}_i^{t}$ that has multiple paths to agent $i$, unless it is among her neighbors (agent $k$ is directly observed by agent $i$).

\begin{proposition}[Graphical Condition for Transparency]\label{PROP:SUFFICIENTCONDITION} 
 Agent $i$ will hold a transparent (efficient) Bayesian posterior belief at time $t$ if for any $k\in \bar{\mathcal{N}}_i^{t}\setminus\mathcal{N}_i$ there is a unique path from $k$ to $i$. 
\end{proposition}

 The graphical condition that is proposed above is only sufficient. For example, structures (\subref{fig:Transparent2a}) and (\subref{fig:Transparent3a}) in Example \ref{example:trasnparency} violate this condition, despite both being transparent. We present the proof of Proposition \ref{PROP:SUFFICIENTCONDITION} in Appendix \ref{App:graphTransparencyProof}. We provide a constructive proof by showing how to compute the Bayesian posterior by aggregating the changes (innovations) in the updated beliefs of neighbors and using the information about beliefs of agents with multiple paths, to correct for redundancies. Accordingly, for structures that satisfy the sufficient condition for transparency, we obtain a simple (and efficient) algorithm for updating beliefs by setting the total innovation at every step equal to the sum of the most recent  innovations observed at each of the neighbors, correcting for those neighbors who are being double-counted. We define innovations as the change in the observed $\log$-belief ratio of agents between two consecutive steps: $\hat{\boldsymbol{\phi}}_{i,t}:= {\boldsymbol{\phi}}_{i,t} - {\boldsymbol{\phi}}_{i,t-1}$,  and initialize them with $\hat{\boldsymbol{\phi}}_{i,0}:= {\boldsymbol{\phi}}_{i,0}  = \boldsymbol{\lambda}_i$.    

\vspace{10pt}

\noindent
 \begingroup\fboxsep=10pt
 \begin{center}
\fbox{
\begin{minipage}{0.9\textwidth}
\parbox{\textwidth}{{\bf Algorithm 3: CORRECTED-INNOVATIONS $(i,t)$}\vspace{0.0cm}\\ 
    {\bf Input:} Graph $\mathcal{G}$ satisfying Proposition \ref{PROP:SUFFICIENTCONDITION}, ${\boldsymbol{\phi}}_{i,t}$, and $\hat{\boldsymbol{\phi}}_{j,t}, j \in \mathcal{N}_i$. \vspace{0.0cm}\\
    {\bf Output:} Posterior $\log$-belief ratio ${\boldsymbol{\phi}}_{i,t+1}$ \vspace{0.05cm}
\begin{enumerate}[leftmargin = *]
    \item AGGREGATE:  $\hat{\boldsymbol{\phi}}_{i,t+1} = \sum\limits_{j \in \mathcal{N}_i} [\hat{\boldsymbol{\phi}}_{j,t} - \sum\limits_{k\in\mathcal{N}_i\cap\mathcal{N}_j^{t}}{\boldsymbol{\phi}}_{k,0} ]$,
    \item UPDATE: ${\boldsymbol{\phi}}_{i,t+1} = {\boldsymbol{\phi}}_{i,t} + \hat{\boldsymbol{\phi}}_{i,t+1}$.
\end{enumerate}
}
\end{minipage}}
\end{center}
\endgroup
\vspace{10pt} 

Note that the transitive networks introduced in Subsection \ref{sec:POSETs}, by definition, satisfy the sufficient condition of Proposition \ref{PROP:SUFFICIENTCONDITION}. Our next corollary summarizes this observation.

\begin{corollary}[Transitivity is sufficient for transparency]\label{COR:TRANSITIVITYANDTRANSPARENCY} All transitive networks are transparent. 
\end{corollary}

Complete graphs are transitive, and therefore, transparent. Directed paths and rooted trees are other classes where Bayesian belief exchange is efficient, since they satisfy the sufficient condition of Proposition \ref{PROP:SUFFICIENTCONDITION}. These special cases are explained next.   

\begin{example}[Complete graphs, directed paths, and rooted trees]\label{example:trees} Complete graphs are a special case where every agent gets to know about the likelihoods of the private signal of all other agents at time one. Subsequently, every agent in a complete graph holds an efficient belief at time two.  Directed paths and rooted (directed) trees are other classes of transparent structures, which satisfy the sufficient structural condition of Proposition \ref{PROP:SUFFICIENTCONDITION}. Indeed, in case of a rooted tree for any agent $k$ that is indirectly observed by agent $i$, there is a unique path connecting $k$ to $i$. As such the correction terms for the sum of innovations in Algorithm 3 is always zero. Hence, for rooted trees we have $\hat{\boldsymbol{\phi}}_{i,t+1} = \sum_{j \in \mathcal{N}_i}\hat{\boldsymbol{\phi}}_{j,t}$: the innovation at each step is equal to the total innovations observed in all the neighbors.
\end{example}

\subsubsection{Efficient Belief Calculations in Transparent Structures}\label{sec:transparentbeliefcalculations}

Here we describe calculations of a Bayesian agent in a transparent structure. Since the network is transparent to agent $i$, she has access to the following information from the beliefs that she has observed in her neighbors at times $\tau\leq t$, before deciding her belief for time $t+1$:

\begin{itemize}[leftmargin=*]
\item Her own signal $\mathbf{s}_i$ and its $\log$-likelihood ratio ${\boldsymbol\lambda}_i$. 
\item Her observations of the neighboring beliefs: $\{\boldsymbol{\mu}_{j,\tau}:j\in\mathcal{N}_i, \tau\leq t\}$. 
\end{itemize}

Due to transparency, the neighboring beliefs reveal the following information about sums of $\log$-likelihood ratios of private signals of subsets of other agents in the network: $ \sum_{k\in\bar{\mathcal{N}}_j^{\tau}}{\boldsymbol\lambda}_k =  {\boldsymbol\phi}_{i,\tau}, \mbox{ for all } \tau \leq t, \mbox{ and any } j \in \mathcal{N}_i$. To decide her belief, agent $i$ constructs the following system of linear equations in $\mbox{card}\left(\bar{\mathcal{N}}_{t+1}\right)+1$ unknowns: $\{{\boldsymbol\lambda}_j:j\in\bar{\mathcal{N}}_{t+1}$, and ${\boldsymbol\phi}^{\star}\}$,  where  ${\boldsymbol\phi}^{\star} = \sum_{j\in\bar{\mathcal{N}}_{t+1}} {\boldsymbol\lambda}_{j}$ is the best possible (statistically efficient) belief for agent $i$ at time $t+1$:
\begin{align}
\begin{cases}
\sum_{k\in\bar{\mathcal{N}}_j^{\tau}}{\boldsymbol\lambda}_k =  {\boldsymbol\phi}_{j,\tau}, \mbox{ for all } \tau \leq t, \mbox{ and any } j \in \mathcal{N}_i,\\
\sum_{j\in\bar{\mathcal{N}}_i^{t+1}} {\boldsymbol\lambda}_{j} - {\boldsymbol\phi}^{\star} =0. \label{eq:linearSystemEquations}
\end{cases}
\end{align} 

Note that \eqref{eq:linearSystemEquations} lists all the information available to agent $i$ when forming her belief in a transparent structure. Hence, transparency is in fact a statement about the linear system of equations in \eqref{eq:linearSystemEquations}: In transparent structures ${\boldsymbol\phi}^{\star}$ can be determined uniquely by solving the linear system \eqref{eq:linearSystemEquations}. Hence, $ {\boldsymbol\phi}_{i,t+1} = {\boldsymbol\phi}^{\star}$, is not only statistically efficient but also computationally efficient. For a transparent structure the complexity of determining  the Bayesian posterior belief at time $t+1$ is the same as the complexity of performing Gauss-Jordan steps which is  $O(n^3)$ for solving the $t \, . \, \mbox{card}(\mathcal{N}_i)$ equations in $\mbox{card}(\bar{\mathcal{N}}_{i}^{t+1})$ unknowns. Note that here we make no attempts to optimize these computations beyond the fact that their growth is polynomial in $n$.

\begin{corollary}[Efficient Computation of Transparent Beliefs]\label{COR:POLYNOMIALTRANSPARENTBELIEFCOMP}
Consider the revealed belief model of opinion exchange in transparent structures. There is an algorithm that runs in polynomial-time and computes the Bayesian posteriors in transparent structures. 
\end{corollary} In general non-transparent cases, the neighboring beliefs are highly non-linear functions of the $\log$-likelihood ratios --- see e.g. \eqref{eq:highlynonlinerbeliefs}, and the above forward reasoning approach can no longer be applied. Indeed, when transparency is violated then beliefs represent what signal profiles agents regard as possible rather than what they know about the $\log$-likelihood ratios of signals of others whom they have directly or indirectly observed.  In particular, the agent cannot use the reported beliefs of the neighbors directly to make inferences about the original causes of those reports which are the private signals. Instead, to keep track of the possible signal profiles that are consistent with her observations the agent employs a version of the IEIS algorithm of Subsection \ref{sec:generalStrcutures} that is tailored to the case of revealed beliefs.

\section{Conclusions, Open Problems, and Future Directions}\label{sec:conclusions}

We proved hardness results for computing Bayesian actions and approximating posterior beliefs in a model of decision making in groups (Theorems \ref{THM:NPHARDNESSACTIONS} and \ref{THM:NPHARDNESSBELIEFS}). 
We also discussed a few generalizations and limitations of those results.
We further augmented these hardness results by offering special cases where Bayesian calculations simplify and efficient computation of Bayesian actions and posterior beliefs is possible (transitive and transparent networks). 

A potentially challenging research direction is to develop a satisfactory theory of rational information exchange in light of computational constraints. It would be
interesting to reconcile more fully our negative results with the more positive picture
presented by \cite{aaronson2005complexity}.
Less ambitiously, a more exact characterization of computational hardness
for different network and utility structures is certainly possible.
Development of an average-case complexity result would be particularly interesting and relevant.

Another major direction is to investigate other configurations and structures for which the computation of Bayesian actions is achievable in polynomial-time, in particular, to develop tight conditions on the network structure that result in necessary and sufficient conditions for transparency. It is also of interest to know the quality of information aggregation; i.e. under what conditions on the signal structure and network topology, Bayesian actions coincide with the best action given the aggregate information of all agents.

\begin{APPENDICES}

\cleardoublepage
\phantomsection
\section{Proof of Theorem \ref{THM:NPHARDNESSACTIONS} (VERTEX-COVER Reduction)}\label{App:vertexConverProof}

Our reduction is from hardness of approximation for the vertex cover problem.  

\begin{definition}[Vertex Cover of a Graph] Given a graph $\hat{\mathcal{G}}_{m,n} =(\hat{\mathcal{V}},\hat{\mathcal{E}})$, with $|\hat{\mathcal{E}}|= m$ edge and $|\hat{\mathcal{V}}| = n$ vertices, a vertex cover $\hat{\Sigma}$ is a subset of vertices such that every edge of $\hat{\mathcal{G}}_{m,n}$ is incident to at least one vertex in
$\hat{\Sigma}$. Let $\hat{\Xi}$ denote the set of all vertex covers of $\hat{\mathcal{G}}_{m,n}$.
\end{definition}

\begin{theorem}[Hardness of approximation of VERTEX-COVER,
  \cite{KMS18}]\label{prob:vertex_cover}
  For every $\varepsilon > 0$,
  given a simple graph $\hat{\mathcal{G}}_{m,n}$ with $n$ vertices and $m$ edges,
  it is NP-hard to distinguish between:
\begin{itemize}
\item YES case: there exists a vertex cover $\hat{\Sigma}$ of size
  $|\hat{\Sigma}| \leq 0.85 n$.
\item NO case: each vertex cover $\hat{\Sigma}$ has size
  $|\hat{\Sigma}| > 0.999n$.
\end{itemize}
\end{theorem}

Theorem~\ref{prob:vertex_cover} follows from recent works on the two-to-two
conjecture culminating in \cite{KMS18}. For completeness, we note that the constants can be improved to
$\sqrt{2}/2+\epsilon$ in the YES case and $1-\epsilon$ in the NO case.

We now restate Theorem~\ref{THM:NPHARDNESSACTIONS} more formally:
\begin{theorem}
  There exists a polynomial-time reduction that maps a graph
  $\hat{\mathcal{G}}_{m,n}$ onto an instance of GROUP-DECISION in the binary
  action model where:
  \begin{itemize}
  \item There are $n+m+2$ agents and the time is set to $t=2$.
  \item For every agent $j$, her signal structure consists of
    efficiently computable numbers that satisfy the following:
    \begin{align}\label{eq:02}
      \exp(-O(n)) < \mathbb{P}_{j,\theta}(1) < 1-\exp(-O(n))\;.
    \end{align}
  \end{itemize}
  Furthermore, letting $i$ be the agent specified in the reduction:
  \begin{itemize}
  \item If $\hat{\mathcal{G}}_{m,n}$ has a vertex cover of size at most
    $0.85n$, then the belief of $i$ at time two satisfies
    $\boldsymbol{\mu}_{i,2}(1) < \exp(-\Omega(n))$.
  \item If all vertex covers of $\hat{\mathcal{G}}_{m,n}$ have size
    at least $0.999n$, then the belief of $i$ satisfies
    $\boldsymbol{\mu}_{i,2}(0) < \exp(-\Omega(n))$.
  \end{itemize}
\end{theorem}

Consider a graph input to the vertex cover problem $\hat{\mathcal{G}}_{m,n}$ with $m$ edges and $n$ vertices.
We encode the structure of $\hat{\mathcal{G}}_{m,n}$ by a two layer network, with $n$ vertex agents $\tau_1, \ldots, \tau_n$ and $m$ edge agents $\varepsilon_1, \ldots, \varepsilon_m$. Each edge agent observes two vertex agents corresponding to its incident vertices in $\hat{\mathcal{G}}_{m,n}$ (see Figure~\ref{fig:VERTEXCOVER}).  Each vertex agent $\tau$ receives a private binary signal $\mathbf{s}_{\tau}$ such that:  
\begin{align}
  &\mathbb{P}_{\tau, 1}(1) = \mathbb{P}_{1}\{\mathbf{s}_{\tau}=1\} =  0.4 =:\overline{p}, \\ 
&\mathbb{P}_{\tau,0}(1)= 0.3 =: \underline{p}, 
\end{align}
where we use the notation $\mathbb{P}_{\theta}\{\cdots\}$ to denote probability
of an event conditioned on the value of the state ($\theta$).

The network also contains two more agents that we call $i$ and $\kappa$.
Agent $\kappa$ does not observe any other agents, while agent $i$
observes $\kappa$ and all edge agents.
We analyze the decision problem of agent $i$ at time $t=2$.
We assume that agent $i$ and the edge agents $\varepsilon_1, \ldots,\varepsilon_m$
receive non-informative private signals.
The signal structure of agent $\kappa$ will be specified later.
We give the observation history of agent $i$ as follows:
All edge agents claim $\mathbf{a}_{\varepsilon_j,1} = 1$ and $\kappa$ claims
$\mathbf{a}_{\kappa,0} = 0$. That concludes the description of the reduction.

Clearly, the reduction is computable in polynomial time and the signal
structures satisfy \eqref{eq:02}, except for agent $\kappa$, which we will check soon. In the rest of the proof, we show that graphs with small vertex covers map onto networks where agent $i$ puts a tiny belief on state one, and graphs with only large vertex covers map onto networks where agent $i$ concentrates her belief on state one.

Consider any edge agent $\varepsilon$ and let $\varepsilon^{(1)}$ and $\varepsilon^{(2)}$ be the vertex agents whose actions are observed by $\varepsilon$. Recalling Example~\ref{example:initial_rounds}, we know
that for any vertex agent $\tau$ her $\log$-belief ratio at time zero
is determined by her private signal:
$\boldsymbol{\phi}_{\tau,0}=\boldsymbol{\lambda}_{\tau}$, and consequently
$\mathbf{a}_{\tau,0}=\mathbf{s}_{\tau}$.

Furthermore, by~\eqref{eq:beliefexchange},
the belief $\boldsymbol{\mu}_{\varepsilon,1}$ and $\log$-belief ratio
$\boldsymbol{\phi}_{\varepsilon,1}$ are determined by the neighboring actions (and private signals) $\mathbf{a}_{\varepsilon^{(1)},0}=\mathbf{s}_{\varepsilon^{(1)}}$ and
$\mathbf{a}_{\varepsilon^{(2)},0}=\mathbf{s}_{\varepsilon^{(2)}}$. Clearly,
if $\mathbf{s}_{\varepsilon^{(1)}}=\mathbf{s}_{\varepsilon^{(2)}}$, then
$\varepsilon$ broadcasts a matching action
$\mathbf{a}_{\varepsilon,1} = \mathbf{s}_{\varepsilon^{(1)}}=\mathbf{s}_{\varepsilon^{(2)}}$.
On the other hand, if $\mathbf{s}_{\varepsilon^{(1)}} \neq \mathbf{s}_{\varepsilon^{(2)}}$, then
the belief of $\varepsilon$ is given by:
\begin{align}
\boldsymbol{\mu}_{\varepsilon,1}(1)  = \frac{\bar{p}(1-\bar{p})}{\bar{p}(1-\bar{p}) + \underline{p}(1-\underline{p})} =  \frac{(0.4)(0.6)}{(0.4)(0.6) + (0.3)(0.7)} > \frac{1}{2},
\end{align}
and therefore $\mathbf{a}_{\varepsilon,1} = 1$, whenever
$\mathbf{a}_{\varepsilon^{(1)},0} \neq \mathbf{a}_{\varepsilon^{(2)},0}$.
To sum up, we have:

\begin{fact}\label{edgeagnetactions}
$\mathbf{a}_{\varepsilon,1} =
\mathds{1}\{\mathbf{s}_{\varepsilon^{(1)}} = 1 \mbox{ or } \mathbf{s}_{\varepsilon^{(2)}} = 1 \}$.
\end{fact} 

The following observation immediately follows from Fact \ref{edgeagnetactions}
and relates our GROUP-DECISION instance to vertex covers of
$\hat{\mathcal{G}}_{m,n}$:

\begin{fact}\label{setcovercondition}
  Define a random variable $\boldsymbol{\Sigma}$ as
  $\boldsymbol{\Sigma}:=\{\tau\in\hat{\mathcal{V}}:\mathbf{s}_{\tau} =1\}$.
  Then, $\boldsymbol{\Sigma}$ is a vertex cover of graph $\hat{\mathcal{G}}_{m,n}=(\hat{\mathcal{V}},\hat{\mathcal{E}})$ if, and only if, $\mathbf{a}_{\varepsilon,1} = 1$ for all $\varepsilon \in \hat{\mathcal{E}}$.
\end{fact}

Recall that we are interested in the decision problem of agent $i$ at time two,
given that she has observed $\mathbf{a}_{\varepsilon,1} = 1$ for all $\varepsilon \in \hat{\mathcal{E}}$, i.e., she has learned that the private signals of vertex agents
form a vertex cover of $\hat{\mathcal{G}}_{m,n}$. Given a particular vertex cover $\hat{\Sigma}$, let us denote its size
by $|\hat{\Sigma}| = \alpha n$ for some
$\alpha=\alpha(\hat{\Sigma})
\in \{\frac{1}{n},\frac{2}{n},\ldots,\frac{n-1}{n},1\}$. Then, we can write
\begin{align}
  \mathbb{P}_{1}\{\boldsymbol{\Sigma}=\hat{\Sigma}\}  
  = \left( \overline{p}^{\alpha} (1-\overline{p})^{(1-\alpha)} \right)^n =:
  \overline{q}(\alpha) ^n,
  \label{eq:03}\\
  \mathbb{P}_{0}\{\boldsymbol{\Sigma}=\hat{\Sigma}\}
  = \left( \underline{p}^{\alpha} (1-\underline{p})^{(1-\alpha)} \right)^n =:  \underline{q}(\alpha) ^n,
  \label{eq:04}
\end{align} where $ \underline{q}(\alpha) =  \underline{p}^{\alpha} (1-\underline{p})^{(1-\alpha)}$ and $\overline{q}(\alpha) =  \overline{p}^{\alpha} (1-\overline{p})^{(1-\alpha)}$. 

We are now ready to consider the Bayesian posterior belief of agent $i$
at time two. It is more convenient to work with the $\log$-belief ratio
$\boldsymbol{\phi}_{i,2}$:
\begin{align}
  \boldsymbol\phi_{i,2}
  &=\log\left(\frac{\boldsymbol{\mu}_{i,2}(1)}{\boldsymbol{\mu}_{i,2}(0)}\right)
    =\log\left(\frac
    {\mathbb{P}\{\theta=1\mid\mathbf{a}_{\varepsilon,1} = 1\mbox{ for all }
    \varepsilon \in \hat{\mathcal{E}}\mbox{ and }\mathbf{a}_{\kappa,0} = 0\}}
    {\mathbb{P}\{\theta=0\mid\mathbf{a}_{\varepsilon,1} = 1\mbox{ for all }
    \varepsilon \in \hat{\mathcal{E}}\mbox{ and }\mathbf{a}_{\kappa,0} = 0\}}
    \right)\\
  &=\log\left(\frac
    {\mathbb{P}_1\{\boldsymbol\Sigma
    \mbox{ is a vertex cover and }\mathbf{a}_{\kappa,0}=0\}}
    {\mathbb{P}_0\{\boldsymbol\Sigma
    \mbox{ is a vertex cover and }\mathbf{a}_{\kappa,0}=0\}}
    \right)
    =\log\left(\frac
    {\displaystyle\sum_{\hat{\Sigma} \in \hat{\Xi}}
    \mathbb{P}_{1}\{\boldsymbol{\Sigma}=\hat{\Sigma}\}
    \mathbb{P}_{\kappa,1}(0)}
    {\displaystyle\sum_{\hat{\Sigma} \in \hat{\Xi}}
    \mathbb{P}_{0}\{\boldsymbol{\Sigma}=\hat{\Sigma}\}
    \mathbb{P}_{\kappa,0}(0)}
    \right)
  \\
  &=\log\left(\frac
    {\displaystyle\sum_{\hat{\Sigma}\in\hat{\Xi}}
    \overline{q}(\alpha)^n}
    {\displaystyle\sum_{\hat{\Sigma} \in \hat{\Xi}}
    \underline{q}(\alpha)^n}
    \right)+\boldsymbol{\lambda}_{\kappa}(0),
\label{eq:logposterior}\end{align}
where along the way we invoked the uniform prior, Fact~\ref{setcovercondition},
as well as~\eqref{eq:03} and~\eqref{eq:04} --- recall $\hat{\Xi}$ denotes the set of all vertex covers of $\hat{\mathcal{G}}_{m,n}$. We now investigate this Bayesian posterior in the YES and NO cases of
VERTEX-COVER. At this point we can also reveal the signal structure of agent $\kappa$:
Letting $q(\alpha):=\log\left(\overline{q}(\alpha)/\underline{q}(\alpha)\right)$,
we choose it such that $\boldsymbol{\lambda}_{\kappa}(0)
  =-({n}/{2})\left(q(0.998)+q(0.999)\right)$
holds (with an arbitrary value for $\boldsymbol{\lambda}_{\kappa}(1)$).

\subsection{Bayesian Posterior in the NO Case:}\label{sec:vertexNOcase}

If we are in the NO case, then all vertex covers have large size
$|\hat{\Sigma}| > 0.999n$.
Since $q(\alpha)=\log\left(\overline{q}(\alpha)/\underline{q}(\alpha)\right)=\alpha\log(14/9)-\log(7/6)$ is a strictly increasing function of
$\alpha$, the expression
$\left(\overline{q}(\alpha)/\underline{q}(\alpha)\right)^n$ also  increases as
$\alpha$ increases. Therefore,~\eqref{eq:logposterior} can be lower-bounded as follows:
\begin{align}
  \boldsymbol\phi_{i,2}
  & = \log\left(\frac
    {\sum_{\hat{\Sigma} \in \hat{\Xi}} \overline{q}(\alpha)^n}{\sum_{\hat{\Sigma} \in \hat{\Xi}} \underline{q}(\alpha)^n }\right) + \boldsymbol\lambda_{\kappa}(0)
    > \log\left(\frac
    {\sum_{\hat{\Sigma} \in \hat{\Xi}}
    \left(\frac{\overline{q}(0.999)}{\underline{q}(0.999)}\cdot\underline{q}(\alpha)\right)^n}
    {\sum_{\hat{\Sigma} \in \hat{\Xi}}
    \underline{q}(\alpha)^n }\right) + \boldsymbol\lambda_{\kappa}(0)
  \\
  &=nq(0.999)-\frac{n}{2}\left(q(0.998)+q(0.999)\right)
    =\frac{n}{2}\left(q(0.999)-q(0.998)\right),
    \label{eq:16}
\end{align}
establishing $\boldsymbol{\phi}_{i,2}=\Omega(n)$
and $\boldsymbol{\mu}_{i,2}(0) = {1}/{(1+e^{\boldsymbol{\phi}_{i,2}})}
  <\exp(-\Omega(n))$,  so that $\mathbf{a}_{i,2}= 1 $.

\subsection{Bayesian Posterior in the YES Case:}\label{sec:vertexYEScase} 
If we are in the YES case, then there exists a small vertex cover $\Sigma^{\star}$ with $|\Sigma^{\star}| = \alpha^{\star} n \leq 0.85 n$. We will show that the total contribution from all large vertex covers with $|\hat{\Sigma}| \geq 0.998 n$ is dominated by the likelihood of this small vertex cover. To this end, we use the following tail bound for the sum of i.i.d. Bernoulli random variables:
\begin{theorem}\label{thm:chernoff}
  Let $\mathbf{s}_1,\ldots,\mathbf{s}_n\in\{0,1\}$ be i.i.d.~binary
  random variables with $\mathbb{P}\{\mathbf{s}_i=1\}=p$ and let
  $p\le\alpha\le 1$. Then,
  \begin{align}\label{eq:06}
    \Pr\left\{
    \sum_{k=1}^n\mathbf{s}_k\ge\alpha n
    \right\}
    \le\exp\left(-nD_{KL}(\alpha||p)
    \right),
  \end{align}
  where the \emph{Kullback-Leibler divergence} $D_{KL}(\cdot)$ is given
  by:
  \begin{align}\label{eq:05}
    D_{KL}(\alpha||p)=\alpha\log\left(\frac{\alpha}{p}\right)
    +(1-\alpha)\log\left(\frac{1-\alpha}{1-p}\right).
  \end{align}
\end{theorem}
The important feature of formula~\eqref{eq:05} is that
$D_{KL}(\alpha||p)$ goes to $\log(1/p)$ as $\alpha$ goes to $1$.
Hence, for every $\delta > 0$ we can choose $\alpha < 1$ such
that the right-hand side of~\eqref{eq:06} is equal to
$(p+\delta)^n$.

In particular, we can upper-bound the likelihood of large
vertex covers with $\theta=1$ as follows: 
\begin{align}\label{eq:08}
  \mathbb{P}_{1}\left\{|\boldsymbol{\Sigma}| = \sum_{k=1}^{n}\mathbf{s}_k>(0.998)n\right\}
  &\leq \exp(- n D_{KL}(0.998||0.4)) = (0.4061...)^n < 0.41^n,
\end{align}  
On the other hand, since $|\Sigma^{\star}| = \alpha^{\star}  n\leq 0.85 n$, we have: 
\begin{align}\label{eq:09}
  \mathbb{P}_{1}\{\boldsymbol{\Sigma}=\Sigma^{\star}\} = \left(\overline{p}^{\alpha^{\star}} (1-\overline{p} )^{1-\alpha^{\star}}\right)^n  \geq  \left(\overline{p}^{0.85} (1-\overline{p})^{0.15}\right)^n  =
  \overline{q}(0.85)^n = \left(0.425\ldots\right)^n > 0.42^n. 
\end{align}
Therefore, in the YES case, conditioned on $\theta=1$ and
$\boldsymbol{\Sigma}$ being a vertex cover,
the probability of having a large vertex cover ($|\boldsymbol{\Sigma}| > 0.998 n$) is exponentially small. 

We are now ready to upper-bound the $\log$-belief ratio in the YES case. Starting again from \eqref{eq:logposterior}, we get:
\begin{align}
  \boldsymbol\phi_{i,2}
  &=\log\left(\frac{\sum_{\hat{\Sigma} \in \hat{\Xi}} \overline{q}(\alpha)^n}{\sum_{\hat{\Sigma} \in \hat{\Xi}}
    \underline{q}(\alpha)^n }\right) + \boldsymbol\lambda_{\kappa}(0)   
  < \log\left(\frac
    {\sum_{\hat{\Sigma} \in \hat{\Xi}:\alpha \leq 0.998} \overline{q}(\alpha)^n
    +\mathbb{P}_1\left\{|\boldsymbol{\Sigma}|>0.998n\right\}}
    {\sum_{\hat{\Sigma} \in \hat{\Xi}: \alpha \leq 0.998} \underline{q}(\alpha)^n}
    \right) + \boldsymbol\lambda_{\kappa}(0) \\
  &< \log\left(\frac
    {\sum_{\hat{\Sigma} \in \hat{\Xi}:\alpha \leq 0.998}
    \left(\frac{\overline{q}(0.998)}{\underline{q}(0.998)}\cdot\underline{q}(\alpha)\right)^n
    \left(1+\left(\frac{0.41}{0.42}\right)^n\right)}
    {\sum_{\hat{\Sigma} \in \hat{\Xi}: \alpha \leq 0.998} \underline{q}(\alpha)^n}
    \right) + \boldsymbol\lambda_{\kappa}(0)
  \label{eq:07}\\ 
  &= nq(0.998) + \log\left(1+\left(\frac{0.41}{0.42}\right)^n\right) 
    -\frac{n}{2}\left(q(0.999)+q(0.998)\right)
    < -\Omega(n),
\end{align}
where we used~\eqref{eq:08} and~\eqref{eq:09} to establish~\eqref{eq:07}.
This implies that $\boldsymbol{\mu}_{i,2}(1) =
{e^{\boldsymbol{\phi}_{i,2}}}/{(1+e^{\boldsymbol{\phi}_{i,2}})} <\exp(-\Omega(n))$,  and $\mathbf{a}_{i,2} = 0$.

From \ref{sec:vertexNOcase} and \ref{sec:vertexYEScase} we conclude that agent $i$ cannot determine her binary action at time two unless she can solve the NP-hard approximation of the VERTEX-COVER problem.   \QEDB

\cleardoublepage
\phantomsection 
\section{Proof of Theorem \ref{THM:NPHARDNESSBELIEFS} (EXACT-COVER Reduction)}\label{App:exactConverProof}

Our reduction is from a variant of the classical EXACT-COVER problem.
An instance of EXACT-COVER consists of a set of $n$ elements
and a collection of sets over those elements. The computational
problem is to decide if there exists a subcollection that exactly covers
the elements, i.e., each element belongs to exactly one set in the
subcollection. We use a restricted version of EXACT-COVER
known as ``Restricted Exact Cover by Three Sets'' (RXC3).

\begin{problem}[RXC3]\label{prob:3dmatching} Consider a set of $n$ elements $\hat{\mathcal{E}}_n$.
Consider also a set $\hat{\mathcal{T}}_n$ of $n$ subsets of
$\hat{\mathcal{E}}_n$, each of them of size three.
Furthermore, assume  that each element
of $\hat{\mathcal{E}}_n$ belongs to exactly three sets in
$\hat{\mathcal{T}}_n$. The RXC3 problem is to decide if there exists a subset $\hat{\mathcal{T}} \subseteq \hat{\mathcal{T}}_n$
  of size $|\hat{\mathcal{T}}| = n/3$ such that it constitutes an exact cover for $\hat{\mathcal{E}}_n$,
  i.e., $\biguplus_{\tau \in \hat{\mathcal{T}}} \tau = \hat{\mathcal{E}}_n$. We refer to instances with and without such an exact cover as YES and NO cases, respectively.
\end{problem}

Note that we make an implicit assumption that $n$ is divisible by three.
It is known that RXC3 is NP-complete.

\begin{theorem}[Section 3 and Appendix A in \cite{Gon85}]\label{lem:3Dmatching}
  RXC3 is NP-complete.
\end{theorem}

Let $\hat{\mathcal{G}}_{n} := (\hat{\mathcal{E}}_n,\hat{\mathcal{T}}_n)$ be
an instance of RXC3.
We encode the structure of $\hat{\mathcal{G}}_{n}$ by a two layer network (cf.~Figure~\ref{fig:EXACTCOVER}).
The first layer is comprised of $n$ agents labeled by the subsets $\tau \in \hat{\mathcal{T}}_n$. The second layer is comprised of $n$ agents labeled by the elements $\varepsilon \in \hat{\mathcal{E}}_n$. Each element agent $\varepsilon$ observes the beliefs of three subset agents,
corresponding to the subsets that contain it.
We denote these three subset agents by
$\varepsilon^{(1)}$, $\varepsilon^{(2)}$, and $\varepsilon^{(3)}$.

The element agents receive non-informative private signals.
Each subset agent $\tau \in \hat{\mathcal{T}}_n$ receives a binary private signal,
$\mathbf{s}_{\tau}\in\{0,1\}$ with the following signal structure:
\begin{align}
\mathbb{P}_{\tau,1}(1) = 1/2 =:\overline{p},\qquad 
\mathbb{P}_{\tau,0}(1) = 1/3 =: \underline{p}. 
\end{align}
Recall our $\log$-likelihood ratio notation from Subsection~\ref{sec:loglike}. Let us define the respective $\log$-likelihood ratios of the zero and one signals as follows:
\begin{align}
  \ell_{1} :=
  \log\left({\overline{p}}/{\underline{p}}\right) = \log(3/2),\qquad
  \ell_{0} :=
  \log\left({(1-\overline{p})}/{(1 -\underline{p})}\right) = \log(3/4). \label{eq:loglikelhiidsignsRXC3}
\end{align}
Under the above definitions, for each subset agent $\tau \in \mathcal{T}_n$ we have: $\boldsymbol\lambda_{\tau} = \mathbf{s}_{\tau}\left(\ell_1 - \ell_0\right) + \ell_0$.

The network contains two more agents called $\kappa$ and $i$.
Agent $\kappa$ is observed by all element agents.  She receives a binary private signal $\mathbf{s}_{\kappa}$
with the following signal structure:
\begin{align}
\overline{p}^{\star} &:= \mathbb{P}_{\kappa,1}(1),  \quad 
\ell^{\star}_{1} := \log\left({\overline{p}^{\star}}/{\underline{p}^{\star}}\right) , \\
\underline{p}^{\star} &:= \mathbb{P}_{\kappa,0}(1),  \quad 
  \ell^{\star}_{0} := \log\left({(1- \overline{p}^{\star})}/{(1 - \underline{p}^{\star})}\right). 
\end{align}
We choose the signal structures such that $\ell^{\star}_{1} - \ell^{\star}_{0} = 2(\ell_1 - \ell_0)$. For concreteness, let $\ell^\star_0:=2\ell_0$ and
$\ell_1^\star:=2\ell_1$.

Finally, agent $i$ does not receive a private signal but observes all element agents (see Figure \ref{fig:EXACTCOVER}). We are interested in the decision problem of agent $i$ at 
time two with the following observation history: Every element agent $\varepsilon\in\hat{\mathcal{E}}_n$ reports the same $\log$-belief ratio at time one:
\begin{align}\label{eq:10}
  \boldsymbol\phi_{\varepsilon,1} = 3\ell_1 + 2\ell_{0}.
\end{align}

Consider the belief of an element agent $\varepsilon$ at time one, given her observations of the subset agents $\varepsilon^{(1)}, \varepsilon^{(2)}, \varepsilon^{(3)}$ and agent $\kappa$. By~\eqref{eq:beliefexchange} and using $\ell_b^\star=2\ell_b$, $b\in\{0,1\}$, we can compute the $\log$-belief ratio of $\varepsilon$ at time one:
\begin{align}
\boldsymbol\phi_{\varepsilon,1} &=  \boldsymbol\lambda_{\kappa} + \sum_{j = 1}^{3}\boldsymbol\lambda_{\varepsilon^{(j)}}  = \mathbf{s}_{\kappa}(\ell^{\star}_{1} - \ell^{\star}_{0}) + \ell^{\star}_{0} + 3\ell_0  + \sum_{j=1}^{3} \mathbf{s}_{\varepsilon^{(j)}}(\ell_{1} - \ell_{0}) \\ 
& = (\ell_1 - \ell_0)(2\mathbf{s}_\kappa + \sum_{j=1}^{3}\mathbf{s}_{{\varepsilon}^{(j)}}) + 5\ell_0 . 
\end{align}

From her observations at time one given by~\eqref{eq:10}, agent $i$ learns that the signals in the neighborhood of each
element agent satisfy the following:
\begin{align}
  2\mathbf{s}_\kappa + \sum_{j=1}^{3}\mathbf{s}_{{\varepsilon}^{(j)}}  = 3.
  \label{eq:systemexactCover}
\end{align}
We denote the set of all signal profiles that satisfy \eqref{eq:systemexactCover} by: 
\begin{align}\label{eq:solutionset}
\Sigma = \left\{(s_{\kappa}, s_{\tau_1}, \ldots , s_{\tau_n}) \in \{0,1\}^{n+1}: 2{s}_\kappa + \sum_{j=1}^{3}{s}_{{\varepsilon}^{(j)}}  = 3, \mbox{ for all } \varepsilon\in\hat{\mathcal{E}}_n\right\}.
\end{align}


Consequently, the $\log$-belief ratio of agent $i$ at time two is given by:
\begin{align}\label{eq:logposteriorratio}
  \boldsymbol\phi_{i,2}
  & =
    \log\left(\frac
    {\displaystyle\sum_{(s_{\kappa}, s_{\tau_1}, \ldots , s_{\tau_n}) \in \Sigma} {(\overline{p}^{\star})}^{s_{\kappa}}{(1-\overline{p}^{\star})}^{1 -s_{\kappa}} {(\overline{p})}^{\sum_{j=1}^{n}s_{\tau_j}}{(1-\overline{p})}^{n -\sum_{j=1}^{n}s_{\tau_j}}} {\displaystyle\sum_{(s_{\kappa}, s_{\tau_1}, \ldots , s_{\tau_n}) \in \Sigma} {(\underline{p}^{\star})}^{s_{\kappa}}{(1-\underline{p}^{\star})}^{1 -s_{\kappa}} {(\underline{p})}^{\sum_{j=1}^{n}s_{\tau_j}}{(1-\underline{p})}^{n -\sum_{j=1}^{n}s_{\tau_j}}}\right).
\end{align}

We now proceed to characterize the solution set $\Sigma$, which determines the posterior ratio per \eqref{eq:logposteriorratio}. One possibility is to set $s_\kappa = 0$, then \eqref{eq:systemexactCover} implies that $s_{{\varepsilon}^{(j)}} = 1$ for all $\varepsilon$ and $j =1 ,2,3$. This is equivalent to having $s_{\tau} =1$ for all subset agents $\tau$. Therefore, $(0,1,1,\ldots,1) \in \Sigma$.

On the other hand, if $s_\kappa = 1$, then \eqref{eq:systemexactCover} implies that:
\begin{align}
  \sum_{j=1}^{3}\mathbf{s}_{{\varepsilon}^{(j)}}  = 1, \mbox{ for every agent }
  \varepsilon\in\hat{\mathcal{E}}_n. \label{eq:systemexactCover2}
\end{align}
In other words, the signal profiles of the subset agents
$({s}_{\tau_1}, {s}_{\tau_2}, \ldots, {s}_{\tau_n})$ specify an exact set-cover of $\hat{\mathcal{E}}_n$. We now investigate the Bayesian posterior of agent $i$ depending on the existence of
an exact set cover.

\subsection{Bayesian Posterior in the NO Case:}\label{sec:exactNOcase}

If we are in a NO case of the RXC3 problem, then the instance $\hat{\mathcal{G}}_{n} = (\hat{\mathcal{E}}_n,\hat{\mathcal{T}}_n)$ does not have an exact set cover.
Therefore,
the solution set $\Sigma$ is a singleton $\Sigma=\{(0,1,\ldots,1)\}$
and~\eqref{eq:logposteriorratio} becomes:
\begin{align}\label{eq:logposteriorratioNOcase}
  \boldsymbol\phi_{i,2}
  & = \log\left(\frac{(1-\overline{p}^{\star}) {(\overline{p})}^{n}} {{(1-\underline{p}^{\star})} {(\underline{p})}^{n}}\right) = 2\ell_0 + n \ell_1
    =2\ell_0 + n\log(3/2) > \Omega(n),
\end{align}
and consequently
$\boldsymbol{\mu}_{i,2}(0)=
{1}/{(1+e^{\boldsymbol{\phi}_{i,2}})} < \exp(-\Omega(n))$.

\subsection{Bayesian Posterior in the YES Case:}\label{sec:exactYEScase}

If we are in a YES case of the RXC3 problem, then there exists
at least one exact cover of $\hat{\mathcal{E}}_n$ consisting of
$n/3$ sets from $\hat{\mathcal{T}}_n$.
Let $\bar{s}' = (s_{\kappa},s_{\tau_1},\ldots,s_{\tau_n}) \in \Sigma$ be a signal configuration corresponding to such an exact cover. Let us also denote the corresponding random profile of private signals by $\bar{\mathbf{s}} = (\mathbf{s}_{\kappa},\mathbf{s}_{\tau_1},\ldots,\mathbf{s}_{\tau_n})$.  The contribution of $\bar{s}$ to the Bayesian posterior of agent $i$ can be calculated as: 
\begin{align}
  &\mathbb{P}_1\{\bar{\mathbf{s}}=\bar{s}'\} = (\overline{p}^{\star}) (\overline{p})^{n/3}{(1-\overline{p})}^{2n/3}= (\overline{p}^{\star}) \overline{q}^{n}, \\
  & \mathbb{P}_0\{\bar{\mathbf{s}}=\bar{s}'\} = (\underline{p}^{\star}) (\underline{p})^{n/3}{(1-\underline{p})}^{2n/3} = (\underline{p}^{\star}) \underline{q}^{n},
\end{align} where 
$\overline{q} := (\overline{p})^{1/3}{(1-\overline{p})}^{2/3} = 1/2$ and
$\underline{q} := (\underline{p})^{1/3}{(1-\underline{p})}^{2/3} \approx 0.529134$.

Let $\hat{N} \geq 1$ be the number of exact covers of $\hat{\mathcal{G}}_n$.
Then, we can use $\overline{p}=\overline{q}$ to compute:
\begin{align}
  \boldsymbol\phi_{i,2} = \log\left(\frac{\hat{N} \cdot \overline{p}^{\star} \cdot \overline{q}^{n} + (1- \overline{p}^{\star}) \cdot \overline{p}^{n}}
  {\hat{N} \cdot \underline{p}^{\star} \cdot \underline{q}^{n} + (1- \underline{p}^{\star}) \cdot\underline{p}^{n}}\right)
  <\log\left(\frac
  {O(\hat{N}\cdot\overline{p}^\star\cdot\overline{q}^n)}
  {\hat{N}\cdot\underline{p}^\star\cdot\underline{q}^n}
  \right)
  \le n\hat{\ell} + O(1)\;
  \le -\Omega(n),
\end{align}
where $\hat{\ell}:= \log(\overline{q}/\underline{q}) < 0$.
Consequently,
$\boldsymbol{\mu}_{i,2}(1)=
{\exp(\boldsymbol{\phi}_{i,2})}/
(1+\exp(\boldsymbol{\phi}_{i,2})) < \exp(-\Omega(n))$.

All in all, from \ref{sec:exactNOcase} and \ref{sec:exactYEScase} we conclude that agent $i$ cannot determine whether her Bayesian posterior concentrates on state zero or state one unless she can solve the NP-hard RXC3 EXACT-COVER variant. \QEDB

\section{I.I.D.~Signals}\label{sec:app:infomrativesignals}

Following Subsection \ref{sec:asymmtericlikelihoods}, we explain how to modify our two reductions (VERTEX-COVER and EXACT-COVER), in Appendices \ref{App:vertexConverProof}  and \ref{App:exactConverProof}, to work with i.i.d. private signals for all agents. 

\subsection{VERTEX-COVER}

Recall that we need to modify our construction such that all 
agents have signal structures of vertex agents, with
$\overline{p}=0.4$, $\underline{p}=0.3$ and respective $\log$-likelihood ratios
$\ell_1=\log(\overline{p}/\underline{p})=\log(4/3)$, and
$\ell_0=\log((1-\overline{p})/(1-\underline{p}))=\log(6/7)$.

The reduction
relies on an auxiliary agent $\kappa$ whose signal $\log$-likelihood ratio is given by
$\boldsymbol\lambda_{\kappa}(0) = -cn$ for some constant $c>0$. Since agent $\kappa$ is directly observed by agent $i$, all we need to do
is to replace $\kappa$ with a number of i.i.d.~agents with the vertex agent
signal structure providing a similar total contribution to $\log$-belief ratio.
Clearly, this is achieved by taking
$n_i = \lfloor cn/ |\ell_0|\rfloor$ agents, all of them broadcasting action
zero (see Figure~\ref{fig:auxillary_n_i}).

We also need to  explain how to handle agents with non-informative signals,
i.e., edge agents and agent $i$. We will leave all those agents in place
and equip them with vertex agent signal structure. We will also indicate in the observation
history that their actions at time zero (and therefore also private signals)
were all ones: $\mathbf{a}_{i,0} = \mathbf{a}_{\varepsilon_1,0} = \ldots = \mathbf{a}_{\varepsilon_m,0} = 1$.

We next add $m$ auxiliary agents $\kappa_1$, $\ldots$, $\kappa_m$
such that each $\kappa_j$ is observed by its corresponding edge agent $\varepsilon_j$, as well as by agent $i$. Again, we let each $\kappa_j$ have the vertex agent signal structure and we indicate that $\mathbf{s}_{\kappa_j}=\mathbf{a}_{\kappa_j,0}=0$. We next verify that Fact \ref{edgeagnetactions} continues to hold. Suppose that an edge agent (called $\varepsilon$) observes opposite actions in her vertex agents (i.e. $\{\mathbf{a}_{\varepsilon^{(1)},0} \, , \, \mathbf{a}_{\varepsilon^{(2)},0}\} = \{0,1\}$). Then the belief of agent $\varepsilon$ at time one aggregates her private signal
(a one signal), the action of her auxiliary agent (a zero signal),
as well as two opposing signals of her vertex agents.
The resulting belief of $\varepsilon$ at time one puts more weight on state one:
\begin{align}
\boldsymbol{\mu}_{\varepsilon,1}(1)  = \frac{\bar{p}^2(1-\bar{p})^2}{\bar{p}^2(1-\bar{p})^2 + \underline{p}^2(1-\underline{p})^2} =  \frac{(0.4)^2(0.6)^2}{(0.4)^2(0.6)^2 + (0.3)^2(0.7)^2} > \frac{1}{2}.
\end{align}
Similarly, if both vertex agents report zero signals, then we see that aggregating three zero signals and one one signal results in a belief that puts less weight on state one:
\begin{align}
\boldsymbol{\mu}_{\varepsilon,1}(1)  = \frac{\bar{p}(1-\bar{p})^3}{\bar{p}(1-\bar{p})^3 + \underline{p}(1-\underline{p})^3} =  \frac{(0.4)(0.6)^3}{(0.4)(0.6)^3 + (0.3)(0.7)^3} < \frac{1}{2}.
\end{align}
Therefore, Fact~\ref{edgeagnetactions} still holds.

The remaining steps of the reduction carry through as before, except that
we need to account for the effect of the new signals of agents $\varepsilon_j$ and
$\kappa_j$, as well as agent $i$'s own private signal. Since these signals amount to $m+1$ ones and $m$ zeros,
their total effect in terms of $\log$-likelihood ratio is equal to $\ell_1 + m(\ell_1+\ell_0)>0$. We can cancel out this net effect  asymptotically by inclusion of $n_i = \left\lfloor {(\ell_1 + m(\ell_1 - \ell_0))}/{|\ell_0|}\right\rfloor$ additional agents that are observed only by agent $i$, each receiving a zero private signal (similar to Figure~\ref{fig:auxillary_n_i}).

\subsection{EXACT-COVER}

In the EXACT-COVER reduction we use an auxiliary agent $\kappa$ whose signal structure is different from those of the subset agents: More precisely, the $\log$-likelihood ratios of the subset agents are $\ell_0$ and $\ell_1$, while for agent $\kappa$
they are $\ell_0^{\star}=2\ell_0$ and $\ell_1^\star=2\ell_1$. Intuitively, we would like to replace agent $\kappa$ with two agents who have the signal structure of the subset agents, and also ensure that the signals of these two agents agree. To achieve this, we use five auxiliary agents $\kappa_1$, $\kappa_2$, $\kappa_3$, $\kappa_4$, and $\kappa_5$ with the signal structure of the subset agents. Suppose every element agent, $\varepsilon_j$, instead of observing agent $\kappa$, observes the two agents $\kappa_1$ and $\kappa_3$. Suppose further that $\kappa_4$ observes $\kappa_1$ and $\kappa_2$; and $\kappa_5$ observes $\kappa_2$ and $\kappa_3$. Finally, let agent $i$, whose decision is NP-hard, observe $\kappa_4$ and $\kappa_5$ (see Figure \ref{fig:auxillary_iid_rx3}).

The private signals of $\kappa_4$ and $\kappa_5$ can be set arbitrarily. Suppose that the belief reports of $\kappa_4$ and $\kappa_5$ at time two implies the following $\log$-belief ratios: $\boldsymbol\phi_{\kappa_4,2} = \boldsymbol\lambda_{\kappa_4} +\ell_0 +\ell_1$ and $\boldsymbol\phi_{\kappa_5,2} = \boldsymbol\lambda_{\kappa_5} +\ell_0 +\ell_1$. From observing $\kappa_4$, agent $i$ learns that the sum of the $\log$-likelihood ratios of the signals of $\kappa_1$ and $\kappa_2$ is $\ell_0+\ell_1$; equivalently, $\kappa_1$ and $\kappa_2$ have received opposite signals. Similarly, from observing $\kappa_5$ agent $i$ learns that $\kappa_2$ and $\kappa_3$ have received opposite signals. Therefore, signals of $\kappa_1$ and $\kappa_3$ must agree. Observing $\kappa_1$ and $\kappa_3$ has the same effect on beliefs of the element agents as observing the single auxiliary agent $\kappa$ with two times the signal strength.
Note that agent $i$ is influenced by what she learned about signals of
$\boldsymbol{\lambda}_{\kappa_2},\boldsymbol{\lambda}_{\kappa_4}$ and
$\boldsymbol{\lambda}_{\kappa_5}$, but this influence is of the order $O(1)$
and therefore does not affect our analysis of Bayesian posteriors
in Appendix~\ref{App:exactConverProof}.

As for the agents without private signals, i.e., agent $i$ and element agents,
the modifications are quite simple. Again, we assume that all these agents
have the same signal structure as the subset agents, and that they report beliefs
at time zero consistent with zero private signals. This introduces a negative shift equal to $(n+1) \ell_0$ in the $\log$-belief ratio of agent $i$,
which can be asymptotically canceled by adding $n_i = \lfloor (n+1)|\ell_0|/\ell_1\rfloor$ more auxiliary agents that
report ones as their signals and are observed only by agent $i$ (similar to Figure~\ref{fig:auxillary_n_i}).
\QEDB

\section{Noisy Actions}\label{app:noisy}

As described in Subsection~\ref{sec:noisy}, let us consider the following noisy variant of the binary
action model: For each computed action $\mathbf{a}_{i,t}$,
the neighboring agents observe the same  action $\mathbf{a}'_{i,t}=\mathbf{a}_{i,t}$ with probability $1-\delta$ and the flipped action
$\mathbf{a}'_{i,t}=1-\mathbf{a}_{i,t}$ with probability $\delta$, for
some $0<\delta<1/2$.

We want to show that Theorem~\ref{THM:NPHARDNESSACTIONS} still holds in this
model, possibly with the constants in the size of the reduction and
the $\exp(-\Omega(n))$ belief approximation factor depending on $\delta$.
To this end, we use the same VERTEX-COVER problem and the
high-level idea as the proof of Theorem~\ref{THM:NPHARDNESSACTIONS}.

Let us start with the general examination of the effect
of noise on an agent $\tau$ that receives a private signal
with signal structure $\mathbb{P}_{\tau,1}(1)=\overline{p}$
and $\mathbb{P}_{\tau,0}(1)=\underline{p}$.
From the perspective of an agent that observes $\tau$, separating
the private signal of $\tau$ and the error in its action does not matter:
All that matters is the likelihood that can be inferred from observing
$\mathbf{a}_{\tau,0}$. The likelihoods of the two possible observations are as follows:
\begin{align}
  \overline{p}'
  &:= \mathbb{P}_1\{\mathbf{a}'_{\tau,0}=1\}
    = \mathbb{P}_1\{\mathbf{a}_{\tau,0}=\mathbf{a}'_{\tau,0}=1
    \mbox{ or } \mathbf{a}_{\tau,0}=1-\mathbf{a}'_{\tau,0}=0
    \}
    = \overline{p}(1-\delta)+(1-\overline{p})\delta,\\
  \underline{p}'
  &:= \mathbb{P}_0\{\mathbf{a}'_{\tau,0}=0\}
    = \underline{p}(1-\delta)+(1-\underline{p})\delta. \label{eq:noisysignalstructures}
\end{align}

From \eqref{eq:noisysignalstructures}, we see that the ``after-noise'' signal
structures are restricted to
$\delta \le \overline{p}',\underline{p}'\le 1-\delta$, rather
than having the full range between $0$ and $1$. Accordingly, we start the VERTEX-COVER reduction by specifying the after-noise signal structures of the vertex agents as follows: $\overline{p}'=1/4+\delta/2$ and
$\underline{p}'=\delta$. It is easy to check that since $\underline{p}' < \overline{p}' < 1/2$ and since
$\underline{p}'(1-\underline{p}') < \overline{p}'(1-\overline{p}')$, an edge agent $\varepsilon$
observing two vertex agents, $\varepsilon^{(1)}$ and $\varepsilon^{(2)}$, still
satisfies the following version of Fact~\ref{edgeagnetactions}:
\begin{align}
  \mathbf{a}_{\varepsilon,1}=1 \mbox{ if, and only if, }
  \mathbf{a}'_{\varepsilon^{(1)},0}=1 \mbox{ or } \mathbf{a}'_{\varepsilon^{(2)}}=1.
  \label{eq:12}
\end{align}
Note that $\mathbf{a}_{\varepsilon,1}$ on the left-hand side is before-noise,
but $\mathbf{a}'_{\varepsilon^{(1)},0}$ and $\mathbf{a}'_{\varepsilon^{(2)},0}$ on the right-hand side, are after-noise. We will now proceed with the analysis of the reduction, encoding vertex
covers of the input graph $\hat{\mathcal{G}}_{n,m}$ in the after-noise actions
$\mathbf{a}'_{\tau,0}$ of vertex agents.

Previously, for each edge of $\hat{\mathcal{G}}_{n,m}$ we placed
an edge agent $\varepsilon$ observing two vertex agents, $\varepsilon^{(1)}$ and $\varepsilon^{(2)}$, corresponding to its incident vertices in $\hat{\mathcal{G}}_{n,m}$.
This time, for each edge we place $k := k(n, m, \delta)$ agents
$\varepsilon(1),\ldots,\varepsilon(k)$, each of them observing the same
two vertex agents and reporting the same noisy actions
$\mathbf{a}'_{\varepsilon(j),1}=1, j=1,\ldots,k$ to agent $i$ (see Figure \ref{fig:auxillary_noisy}).
Since by~\eqref{eq:12} the before-noise actions $\mathbf{a}_{\varepsilon(j),1}$
have all been the same,  after observing $\varepsilon(1),\ldots,\varepsilon(k)$ agent $i$ concludes
that exactly one of the following is true:
\begin{itemize}
\item Nobody among $\varepsilon(1)$, $\ldots$, $\varepsilon(k)$ has flipped her action:
  $\mathbf{a}_{\varepsilon(j),1} = \mathbf{a}'_{\varepsilon(j),1} = 1$, $j=1,\ldots,k$.
\item Everybody in $\varepsilon(1)$, $\ldots$, $\varepsilon(k)$ has flipped her action:
  $\mathbf{a}_{\varepsilon(j),1} = 1-\mathbf{a}'_{\varepsilon(j),1} = 0$, $j=1,\ldots,k$.
\end{itemize}
In the second case, we say that ``an \emph{error} has occurred in edge $\varepsilon$''.
We now proceed to show that for $k$ large enough the probability of an error occurring in some edge is so small that the analysis of the noisy model essentially reduces back to what we did in Appendix~\ref{App:vertexConverProof}.

To this end, consider the $\log$-belief ratio of agent $i$ at time two,
as described in~\eqref{eq:logposterior}, neglecting for the moment
the influence of agent $\kappa$. Let $\boldsymbol{\Sigma}:=\{\tau\in\hat{\mathcal{V}}:\mathbf{a}'_{\tau,0}=1\}$.
Following Appendix~\ref{App:vertexConverProof}, agent $i$ concludes that either $\boldsymbol{\Sigma}$ forms a vertex
cover of $\hat{\mathcal{G}}_{n,m}$ or an error has occurred in at least one edge. Let $\mathcal{E}$ denote the event that at least one error has occurred and let $\lnot\mathcal{E}$ denote its complement. If we want to, for example, upper-bound
$\boldsymbol{\phi}_{i,2}$ in the NO case, we can write:
\begin{align}\label{eq:14}
  \boldsymbol{\phi}_{i,2}
  &=\log\left(\frac{\boldsymbol{\mu}_{i,2}(1)}{\boldsymbol{\mu}_{i,2}(0)}\right)
    \le \log\left(\frac
    {\mathbb{P}_1\{\boldsymbol{\Sigma}\text{ is a vertex cover}\land\lnot\mathcal{E}\}
    +\mathbb{P}_1\{\mathcal{E}\mid\mathbf{h}_{i,2}\}}
    {\mathbb{P}_0\{\boldsymbol{\Sigma}\text{ is a vertex cover}\land\lnot\mathcal{E}\}}
    \right),
\end{align}
where we drop the error probability term ($\mathbb{P}_0\{\mathcal{E}\mid\mathbf{h}_{i,2}\}$) in the denominator to obtain an upper-bound. We now note that, by union bound and other elementary considerations, the
error probability can be bounded by:
\begin{align}
  \mathbb{P}_1\{\mathcal{E}\mid\mathbf{h}_{i,2}\}
  &\le m \delta^k \left(\delta^k+(1-\delta)^k\right)^{m-1}
  = m \left(\frac{\delta^k+(1-\delta)^k}{(1-\delta)^k}\right)^{m-1}
  \left(\frac{\delta}{1-\delta}\right)^k(1-\delta)^{km}\\
  &\le m2^m\left(\frac{\delta}{1-\delta}\right)^k(1-\delta)^{km}
    \le (1/8+\delta/4)^n (1-\delta)^{km} = (\overline{p}'/2)^n  (1-\delta)^{km}\\
  &=(1/2)^n \mathbb{P}_1\{\forall\tau: \mathbf{a}'_{\tau,0}=1\land\lnot\mathcal{E}
    \}
    \le(1/2)^n\mathbb{P}_1\{\boldsymbol{\Sigma}\text{ is a vertex cover}\land\lnot\mathcal{E}\}
    ,\label{eq:15}
\end{align} where in the first line we use the fact that conditioned on observation
history $\mathbf{h}_{i,2}$, for each edge either $0$ or $k$ flips has occurred.
In the second line, we make $m2^m\left({\delta}/({1-\delta})\right)^k
    \le (1/8+\delta/4)^n$ by choosing $k$ to be (polynomially) large enough --- recall ${\delta}/({1-\delta}) < 1$.

Taken together, \eqref{eq:14} and~\eqref{eq:15} imply that, up to a tiny
$\exp(-\Omega(n))$ factor, the value of $\boldsymbol{\phi}_{i,2}$ is almost
the same as that computed in~\eqref{eq:logposterior}. Therefore,
we can use the same computation as in~\eqref{eq:16} to establish
a linear lower-bound on $\boldsymbol{\phi}_{i,2}$. The YES case is handled very
similarly.

Finally, it remains to account for agent $\kappa$. This is done in basically
the same way as in Subsection~\ref{sec:asymmtericlikelihoods}: We 
replace agent $\kappa$ with the strong $\log$-likelihood
 ratio $\boldsymbol{\lambda}_{\kappa}(0) = -cn$ by
$\left\lfloor\frac{cn}{\delta'}\right\rfloor$ agents with the after-noise $\log$-likelihood ratio
$\boldsymbol{\lambda}'_{\kappa_j}(0) = -\delta'$ for appropriately small
$\delta' = \delta'(\delta)$, all reporting action zero at time zero. This concludes the description of our modification
in the noisy model.

\section{Complexity of Bayesian Decisions Using Algorithm 1: IEIS}\label{App:BAYESGROUPCOMP}
Suppose that agent $i$ has reached her $t$-th decision in a general network structure. Given her information at time $t$, for all $\overline{s} = (s_1,\ldots,s_n)\in\mathcal{S}_1\times\ldots\times\mathcal{S}_n$ and any $j \in {\mathcal{N}}_i^{\tau}, \tau= t+1, t, \ldots, 1$ she has to update $\mathcal{I}(j,t-\tau,\overline{s})$ into $\mathcal{I}(j,t+1-\tau,\overline{s})\subset \mathcal{I}(j,t-\tau,\overline{s})$.  If $\tau = t+1$ then agent $j \in {\mathcal{N}}_i^{\tau}$ is being considered for the first time at the $t$-th step and $\mathcal{I}(j,0,\overline{s}) = \{s_j\}\times\prod_{k\neq j}\mathcal{S}_k$ is initialized without any calculations. However if $\tau \leq t$, then $\mathcal{I}(j,t-\tau,\overline{s})$ can be updated into $\mathcal{I}(j,t+1-\tau,\overline{s})\subset \mathcal{I}(j,t-\tau,\overline{s})$ only by verifying the condition $a_{k,t-\tau}(\overline{s}') = a_{k,t-\tau}(\overline{s})$ for every $\overline{s}' \in \mathcal{I}(j,t-\tau,\overline{s})$ and $k\in\mathcal{N}_j$: any $\overline{s}'\in \mathcal{I}(j,t-\tau,\overline{s})$ that violates this condition for some $k\in\mathcal{N}_j$ is eliminated and  $\mathcal{I}(j,t+1-\tau,\overline{s})$ is thus obtained by pruning $\mathcal{I}(j,t-\tau,\overline{s})$. 

Verification of $a_{k,t-\tau}(\overline{s}') = a_{k,t-\tau}(\overline{s})$  involves calculations of $a_{k,t-\tau}(\overline{s}')$ and $a_{k,t-\tau}(\overline{s})$ according to \eqref{eq:ActionsGivenSignalProfiles}. The latter requires the addition of $\mbox{card}(\mathcal{I}(k,t-\tau,\overline{s}))$ product terms $u_k(a_k,{\theta'})$ ${\mathbb{P}}_{{\theta'}}(\overline{s}')$ $\nu({\theta'})$ $ = $ $u_k(a_k,{\theta'})$ ${{\mathbb{P}}_{1,\theta'}}({s}'_1)$ $\ldots$ ${\mathbb{P}}_{n,{\theta'}}({s}'_n)$ $\nu({\theta}')$ for each $\overline{s}' \in \mathcal{I}(k,t-\tau,\overline{s})$, ${\theta}'\in\Theta$, and $a_k \in \mathcal{A}_k$ to evaluate the left hand-side of \eqref{eq:ActionsGivenSignalProfiles}. Hence, we can estimate the total number of additions and multiplications required for calculation of each conditional action $a_{k,t-\tau}(\overline{s})$ as $A \, . \, (n+2) \, . \, m \, . \, \mbox{card}(\mathcal{I}(k,t-\tau,\overline{s}))$, where $m:=\mbox{card}(\Theta)$ and $A = \max_{k\in[n]}\mbox{card}(\mathcal{A}_k)$. Hence the total number of additions and multiplications undertaken by agent $i$ at time $t$ for determining actions $a_{k,t-\tau}(\overline{s})$ can be estimated as follows:
\begin{align}
C_1 := A \, . \, (n+2) \, . \, \mbox{card}(\Theta) \, . \, \sum_{j\in\bar{\mathcal{N}}_i^t}\sum_{k\in\mathcal{N}_j} \mbox{card}(\mathcal{I}(k,t-\mbox{dist}(j,i),\overline{s}))  \leq A \, . \, (n+2) \, . \, n \, . \, M^{n-1} \, . \, m, \label{calculatins}
\end{align} where we upper-bound the cardinality of the union of the higher-order neighborhoods of agent $i$ by the total number of agents: $\mbox{card}(\bar{\mathcal{N}}_i^{t+1}) \leq n$ and use the inclusion relationship $\mathcal{I}(k,t-\mbox{dist}(j,i),\overline{s}) \subset \mathcal{I}(k,0,\overline{s}) = \{s_k\}\times\prod_{j\neq k}S_j$ to upper-bound  $\mbox{card}(\mathcal{I}(k,t-\mbox{dist}(j,i),\overline{s}))$ by $M^{n-1}$ where $M$ is the largest cardinality of  finite signal spaces, $S_j, j\in[n]$ . As the above calculations are performed at every signal profile $\overline{s} \in \mathcal{S}_1\times\ldots\mathcal{S}_n$ the total number of calculations (additions and multiplications) required for the Bayesian decision at time $t$, denoted by $C_t$, can be bounded  as  follows:
\begin{align}
A \, . \, M^{n} \leq  C_t   \leq A \, . \,  (n+2) \, . \, n \, . \, M^{2n-1} \, . \, m, \label{BayesianCalculations}
\end{align} where on the right-hand side, we apply \eqref{calculatins} for each of the $M^n$ signal profiles. In particular, the calculations grow exponential in the number of agents $n$. Once agent $i$ calculates the actions $a_{k,t-\tau}(\overline{s})$ for all $k\in\mathcal{N}_j$  she can then update the possible signal profiles $\mathcal{I}(j,t-\tau,\overline{s})$, following step 2(a)ii  of Algorithm 1, to obtain $\mathcal{I}(j,t+1-\tau,\overline{s})$ for all $j\in \bar{\mathcal{N}}^{t}_i$ and any $\overline{s}\in\mathcal{S}_1\times\ldots\times\mathcal{S}_n$. This in turn enables her to calculate the conditional actions of her neighbors $a_{j,t}(\overline{s})$ at every signal profile and to eliminate any $\overline{s}$ at which the conditional action   $a_{j,t}(\overline{s})$ does not agree with the observed action $\mathbf{a}_{j,t}$, for some $j\in\mathcal{N}_i$. She can thus update her list of possible signal profiles from $\bm{\mathcal{I}}_{i,t}$ to $\bm{\mathcal{I}}_{i,t+1}$ and adopt the corresponding Bayesian belief ${\boldsymbol\mu}_{i,t+1}$ and action $\mathbf{a}_{i,t+1}$. The latter involves an additional $(n+2)mA\cdot\mbox{card}(\bm{\mathcal{I}}_{i,t+1})$ additions and multiplication which are, nonetheless, dominated by the number of calculations required in \eqref{BayesianCalculations} for the simulation of other agents' actions at every signal profile. \QEDB

\section{Computational Complexity of Algorithm 2: IEIS-TRANSITIVE}\label{App:complexityOFposets}

According to (I2), in a transitive structure at time $t$ agent $i$ has access to the list of possible private signals for each of her neighbors: $\bm{\mathcal{S}}_{j,t}$, $j\in\mathcal{N}_i$ given their observations up until that point in time. The possible signal set $\bm{\mathcal{S}}_{j,t}$ for each agent $j\in\mathcal{N}_i$ is calculated based on the actions taken by others and observed by agent $j$ until time $t-1$ together with possible private signals that can explain her history of choices: $\mathbf{a}_{j,0}$,  $\mathbf{a}_{j,1}$, and so on up until her most recent choice which is $\mathbf{a}_{j,t}$. At time $t$, agent $i$ will have access to all the observations of every agent in her neighborhood and can vet their most recent choices $\mathbf{a}_{j,t}$ against their observations to eliminate the incompatible private signals from the possible set $\bm{\mathcal{S}}_{j,t}$ and obtain an updated list of possible signals  $\bm{\mathcal{S}}_{j,t+1}$ for each of her neighbors $j\in\mathcal{N}_i$. This pruning is achieved by calculating $\mathbf{a}_{j,t}(s_j)$ given $\bm{\mathcal{I}}_{j,t}(s_j) = \{s_j\}\times\prod_{k\in\mathcal{N}_j}{\bm{\mathcal{S}}_{j,t}}$ for each $s_j \in {\bm{\mathcal{S}}_{j,t}}$ and removing any incompatible $s_j$ that violates the condition $\mathbf{a}_{j,t} = {\mathbf{a}_{j,t}}(s_j)$; thus obtaining the pruned set $\bm{\mathcal{S}}_{j,t+1}$. The calculation of ${\mathbf{a}}_{j,t}(s_j)$ given $\bm{\mathcal{I}}_{j,t}(s_j) = \{s_j\}\times\prod_{k\in\mathcal{N}_j}{\bm{\mathcal{S}}_{j,t}}$ is performed according to \eqref{eq:ActionsGivenSignalProfiles} but the decomposition of the possible signal profiles based on the relation $\bm{\mathcal{I}}_{j,t}(s_j) = \{s_j\}\times\prod_{k\in\mathcal{N}_j}{\bm{\mathcal{S}}_{j,t}}$ together with the independence of private signals across different agents help reduce  the number of additions and multiplications involved as follows: 

\begin{align}
\mathcal{A}_j(\bm{\mathcal{I}}_{j,t}(s_j)) & = \argmax_{a_j\in\mathcal{A}_j}\sum_{{\theta'}\in \Theta}u_j(a_j,{\theta'})\frac{\sum_{\overline{s}'\in\bm{\mathcal{I}}_{j,t}(s_j)} {\mathbb{P}}_{\theta'}(\overline{s}')\nu({\theta'})}{\sum_{{\theta''}\in\Theta}\sum_{\overline{s}'\in\bm{\mathcal{I}}_{j,t}(s_j)} {\mathbb{P}}_{\theta''}(\overline{s}')\nu({\theta''})} \\ & = \argmax_{a_j\in\mathcal{A}_j}\sum_{{\theta'}\in\Theta}u_j(a_j,{\theta'})\frac{{\mathbb{P}}_{\theta'}(s_j)\prod_{k\in\mathcal{N}_j}\sum_{s_k\in\bm{\mathcal{S}}_{k,t}} {\mathbb{P}}_{\theta'}(s_k)\nu({\theta'})}{\sum_{{\theta''}\in\Theta}{\mathbb{P}}_{\theta''}(s_j)\prod_{k\in\mathcal{N}_j}\sum_{s_k\in\bm{\mathcal{S}}_{k,t}} {\mathbb{P}}_{\theta''}(s_k)\nu({\theta''})}. 
\label{eq:ActionsGivenSignalProfilesDecomposed}
\end{align}

Hence, the calculation of the conditionally feasible action $\mathbf{a}_{j,t}(s_j)$ for each $s_j\in\bm{\mathcal{S}}_{j,t}$ can be achieved through $\mbox{card}(\Theta)$ $A$ $\sum_{k\in\mathcal{N}_j}\mbox{card}(\bm{\mathcal{S}_{k,t}})$ additions and $\mbox{card}(\Theta)\left(\mbox{card}({\mathcal{N}_j})+2\right)A$ multiplications; subsequently, the total number of additions and multiplications required for agent $i$ to update the possible private signals of each of her neighbor can be estimated as follows:
\begin{align}
A\sum_{j\in\mathcal{N}_i}\mbox{card}(\Theta)\mbox{card}(\bm{\mathcal{S}}_{j,t})\left[\sum_{k\in\mathcal{N}_j}\mbox{card}(\bm{\mathcal{S}_{k,t}}) + \mbox{card}(\mathcal{N}_j) + 2\right] \\ \leq An^2M^2m + An^2Mm + 2nMmA, 
\label{eq:POSETCausalCalculations} 
\end{align} where $M$, $n$, $m$ and $A$ are as in \eqref{BayesianCalculations}. After updating her lists for the possible signal profiles of all her neighbors, the agent can refine her list of possible signal profiles $\bm{\mathcal{I}}_{i,t+1} = \{\mathbf{s}_i\}\times\prod_{j\in\mathcal{N}_i}{\bm{\mathcal{S}}_{j,t+1}}$ and determine her belief ${\boldsymbol\mu}_{i,t+1}$ and refined choice $\mathbf{a}_{i,t+1}$. The latter is achieved through an extra $\mbox{card}(\Theta) A  \sum_{j\in\mathcal{N}_i}\mbox{card}(\bm{\mathcal{S}_{j,t+1}})$ additions and $\mbox{card}(\Theta)A\left(\mbox{card}({\mathcal{N}_i})+2\right)$ multiplications, which are dominated by the required calculations in \eqref{eq:POSETCausalCalculations}. Most notably, the computations required of the agent for determining her Bayesian choices in a transitive network increase polynomially in the number of agents $n$, whereas in a general network structure using IEIS these computations increase exponentially fast in the number of agents $n$.  \QEDB

\section{Proof of Proposition \ref{PROP:SUFFICIENTCONDITION} (Graphical Condition for Transparency)}\label{App:graphTransparencyProof}

The proof follows by induction on $t$, i.e. by considering the agents whose information reach agent $i$ for the first time at $t$. The claim is trivially true at time one, since agent $i$ can always infer the likelihoods of the private signals of each of her neighbors by observing their beliefs at time one. Now consider the belief of agent $i$ at time $t$, the induction hypothesis implies that $\boldsymbol{\phi}_{i,t-1} = \sum_{k\in\bar{\mathcal{N}}_i^{t-1}}\boldsymbol{\lambda}_k$, as well as  $\boldsymbol{\phi}_{j,t-1} = \sum_{k\in\bar{\mathcal{N}}_j^{t-1}}\boldsymbol{\lambda}_k$ and $\boldsymbol{\phi}_{j,t-2} = \sum_{k\in\bar{\mathcal{N}}_j^{t-2}}\boldsymbol{\lambda}_k$ for all $j \in \mathcal{N}_i$.
To form her belief at time $t$ (or equivalently its $\log$-belief ratio $\boldsymbol{\phi}_{i,t}$), agent $i$  should consider her most recent information $\{\boldsymbol{\phi}_{j,t-1} = \sum_{k\in\bar{\mathcal{N}}_j^{t-1}}\boldsymbol{\lambda}_k,j\in\mathcal{N}_i\}$ and use that to update her current belief $\boldsymbol{\phi}_{i,t-1} = \sum_{k\in\bar{\mathcal{N}}_i^{t-1}}\boldsymbol{\lambda}_k$. To prove the induction claim, it suffices to show that agent $i$ has enough information to calculate the sum of $\log$-likelihood ratios of all signals in her $t$-radius ego-net, $\bar{\mathcal{N}}_i^{t}$; i.e. to form $\boldsymbol{\phi}_{i,t} = \sum_{k\in\bar{\mathcal{N}}_i^{t}}\boldsymbol{\lambda}_k$. This is the best possible belief that she can hope to achieve at time $t$, and it is the same as her Bayesian posterior, had she direct access to the private signals of all agents in her $t$-radius ego-net. To this end, by using her knowledge of ${\boldsymbol{\phi}}_{j,t-1}$ and ${\boldsymbol{\phi}}_{j,t-2}$ she can form:
\begin{align} 
\hat{\boldsymbol{\phi}}_{j,t-1} = {\boldsymbol{\phi}}_{j,t-1} - {\boldsymbol{\phi}}_{j,t-2} =  \sum_{k\in\mathcal{N}_j^{t-1}}\boldsymbol{\lambda}_k,
\end{align} for all $j\in\mathcal{N}_i$. Since, $\boldsymbol{\phi}_{i,t-1} = \sum_{k\in\bar{\mathcal{N}}_i^{t-1}}\boldsymbol{\lambda}_k$ by the induction hypothesis, the efficient belief $\boldsymbol{\phi}_{i,t} = \sum_{k\in\bar{\mathcal{N}}_i^{t}}\boldsymbol{\lambda}_k$ can be calculated if and only if,  
\begin{align} 
\hat{\boldsymbol{\phi}}_{i,t} = {\boldsymbol{\phi}}_{i,t} - {\boldsymbol{\phi}}_{i,t-1} =  \sum_{k\in{\mathcal{N}}_i^{t}}\boldsymbol{\lambda}_k, \label{eq:innnovationterm}
\end{align} can be computed. In the above formulation $\hat{\boldsymbol{\phi}}_{i,t}$ is an innovation term, representing the information that agent $i$ learns from her most recent observations at time $t$. We now show that under the assumption that any agent with multiple paths to an agent $i$ is directly observed by her, the innovation term in \eqref{eq:innnovationterm} can be constructed from the knowledge of  $\boldsymbol{\phi}_{j,t-1} = \sum_{k\in\bar{\mathcal{N}}_j^{t-1}}\boldsymbol{\lambda}_k$, and $\boldsymbol{\phi}_{j,t-2} = \sum_{k\in\bar{\mathcal{N}}_j^{t-2}}\boldsymbol{\lambda}_k$ for all $j \in \mathcal{N}_i$; indeed, we show that: 
 \begin{align} 
\hat{\boldsymbol{\phi}}_{i,t} = \sum_{j \in \mathcal{N}_i}\left(\hat{\boldsymbol{\phi}}_{j,t-1} - \sum_{k\in\mathcal{N}_i\cap\mathcal{N}^{t-1}_j}{\boldsymbol{\phi}}_{k,0} \right), \mbox{ for all } t > 1. \label{UpdatingTheInnovationTerms}
\end{align}
Consider any $k \in {\mathcal{N}}_i^{t}$, these are all agents which are at distance exactly $t$, $t>1$, from agent $i$, and no closer to her. No such $k \in {\mathcal{N}}_i^{t}$ is a direct neighbor of agent $i$ and the structural assumption therefore implies that there is a unique neighbor of agent $i$, call this unique neighbor $j_k \in \mathcal{N}_i$, satisfying $k \in {\mathcal{N}}_{j_k}^{t-1}$. On the other hand, consider any $j \in \mathcal{N}_i$ and some $k \in \mathcal{N}_{j}^{t-1}$, contributing $\boldsymbol{\lambda}_k$ to $\hat{\boldsymbol{\phi}}_{j,t-1}$. Such an agent $k$ is either a neighbor of $i$ or else at distance exactly $t>1$ from agent $i$ and therefore $k \in \mathcal{N}_{i}^{t}$, and element $j$ would be the unique neighbor $j_k \in\mathcal{N}_i$ satisfying $k \in {\mathcal{N}}_{j_k}^{t-1}$. Subsequently, using the notation $\uplus$ for disjoint unions, we can partition $\mathcal{N}_{i}^{t}$ as follows: 
\begin{align}
\mathcal{N}_{i}^{t} = \uplus_{j\in\mathcal{N}_i}\bar{\mathcal{N}}_{j}^{t-1}\setminus \left(\bar{\mathcal{N}}_{j}^{t-2}\cup\mathcal{N}_i\right),
\end{align} and therefore we can rewrite the left-hand side of \eqref{eq:innnovationterm} as follows:
\begin{align} 
\hat{\boldsymbol{\phi}}_{i,t} & =  \sum_{k\in {\mathcal{N}}_i^{t}}\boldsymbol{\lambda}_k =  \sum_{\substack{k\in\uplus_{j\in\mathcal{N}_i} \\ \bar{\mathcal{N}}_{j}^{t-1}\setminus (\bar{\mathcal{N}}_{j}^{t-2}\cup\mathcal{N}_i)}}\boldsymbol{\lambda}_k =  \sum_{j\in\mathcal{N}_i} \sum_{\substack{k\in\bar{\mathcal{N}}_j^{t-1}\setminus \\ (\bar{\mathcal{N}}_j^{t-2}\cup\mathcal{N}_i)}}\boldsymbol{\lambda}_k  \\  & = \sum_{j\in\mathcal{N}_i} \left(\sum_{k\in{\mathcal{N}}_j^{t-1}}\boldsymbol{\lambda}_k -  \sum_{k\in\mathcal{N}_i \cap {\mathcal{N}}_j^{t-1}}\boldsymbol{\lambda}_k \right) =  \sum_{j \in \mathcal{N}_i}\left(\hat{\boldsymbol{\phi}}_{j,t-1} - \sum_{k\in\mathcal{N}_i\cap \mathcal{N}_j^{t-1}}{\boldsymbol{\phi}}_{k,0} \right), \label{eq:innnovationtermREwrite}
\end{align} as claimed in \eqref{UpdatingTheInnovationTerms}, completing the proof. \QEDB

\end{APPENDICES}

\printendnotes

\ACKNOWLEDGMENT{This work was partially supported by awards ONR N00014-16-1-2227, NSF CCF-1665252 and ARO MURIs W911NF-12-1-0509, W911NF-19-0217 and by a Vannevar Bush Fellowship.}

\bibliographystyle{informs2014} 
\bibliography{BayesRef}

\end{document}